\documentclass[11pt,twoside]{amsart}

\usepackage{enumerate}
\usepackage{graphicx}
 \usepackage{wrapfig}
\usepackage[hyperindex=true,plainpages=false,colorlinks=false,pdfpagelabels]{hyperref}
\DeclareMathOperator{\End}{End} 
\DeclareMathOperator{\GL}{GL} 
 
\DeclareMathOperator{\Span}{Span}

\renewcommand{\Im}{\operatorname{Im}}

\newcommand{\R}{\mathbb{R}}
\newcommand{\C}{\mathbb{C}}

\newcommand{\Z}{\mathbb{Z}}
\renewcommand{\P}{\mathbb{P}}
\renewcommand{\H}{\mathbb{H}}  
\newcommand{\RP}{\mathbb{RP}}
\newcommand{\CP}{\mathbb{CP}}
\newcommand{\HP}{\mathbb{HP}}
\renewcommand{\S}{S}

\newcommand{\eprint}[1]{e--print: \href{http://#1}{\nolinkurl{#1}}}

\theoremstyle{plain}
\newtheorem{The}{Theorem}[section]
\newtheorem{Pro}[The]{Proposition}
\newtheorem{Lem}[The]{Lemma}
\newtheorem{Cor}[The]{Corollary}
\newtheorem*{Cor*}{Corollary}

\theoremstyle{definition}
\newtheorem{Def}[The]{Definition}

\theoremstyle{remark} 
\newtheorem{Rem}[The]{Remark}

\newtheorem*{Rem*}{Remark}

\numberwithin{equation}{section}
\begin{document}

\title[On Discrete Differential Geometry in Twistor Space]{On Discrete Differential Geometry in Twistor Space}

\date{\today}

\author{George Shapiro}

\address{George Shapiro\\
  41 Montague Avenue\\
  Lake Pleasant, MA 01376}
  
\email{geoshapiro@gmail.com}

\begin{abstract}
 In this paper we introduce a discrete integrable system generalizing the discrete (real) cross-ratio system in $S^4$ to complex values of a generalized cross-ratio by considering $S^4$ as a real section of the complex Pl\"ucker quadric, realized as the space of two-spheres in $S^4.$   We develop the geometry of the Pl\"ucker quadric by examining the novel contact properties of two-spheres in $S^4,$ generalizing classical Lie geometry in $S^3.$   Discrete differential geometry aims to develop discrete equivalents of the geometric notions and methods of classical differential geometry.  We define discrete principal contact element nets for the Pl\"ucker quadric and prove several elementary results. Employing a second real real structure, we show that these results generalize previous results by Bobenko and Suris $(2007)$ on discrete differential geometry in the Lie quadric. 
\end{abstract}

\maketitle

\section{Introduction}\label{sec:intro}

The study of special surfaces in three-dimensional space has been a topic of interest going back to the foundations of differential geometry.  The subject encompasses surfaces in Euclidean space  with familiar properties e.g. constant mean curvature, constant gaussian curvature and surfaces in classical affine, M\"obius (conformal), and projective geometries and also surfaces in  unfashionable \cite{unf} geometries  such as Lie geometry \cite{cec}.   This view was summarized at the beginning of the $20^{th}$ century in the compendium of Blaschke \cite{bla}.  Parameterizations of surfaces in space are determined by solutions of partial differential equations.  Then,  integrability conditions on these equations for special surfaces determine completely integrable partial differential equations. This was an essential component of classical results and was rediscovered in the latter part of $20^{th}$ century with the redevelopment of ``soliton'' geometry \cite{sym}.  

It is a characteristic of the smooth theory that surfaces in $3$-space can often be found in families related by fundamental transformations such as the Ribaucour transformation \cite{bur}.  Under the rubric of ``Discrete Differential Geometry,'' recent work has discretized classical surface theory by modelling the special differential equations of surface theory with partial {\em difference} equations \cite{bob1}, \cite{dol6}.  This has required the development of a discrete theory of ``integrability'' \cite{nov}  defined by modelling the permutability relations (e.g. Bianchi permutability) between transformations of smooth solutions of integrable systems as a ``consistency'' condition in the discrete case \cite{bob1}.  Then, one is lead to the viewpoint  that the geometry of the fundamental transformations and the geometry of the {\em discrete} surface are the same \cite{dol2}.   


Given a parameterization of a surface in space with coordinates $(u,v),$ a system of  coordinate curves defined by $u = $ const. and $v =$ const. determines a ``net'' on the surface.  The classical approach to special surfaces starts with the geometry of these coordinate curves. In particular,  ``conjugate nets,'' are parameterizations of surfaces in space defined by the property that embedded tangent vectors at each point on a coordinate curve stay in the tangent plane  under motions of the curve in the transverse coordinate direction.  This is a concept of projective differential geometry and in an affine chart is equiavalent to the property that the mixed partial derivative of the parameterization is, at each point, spanned by tangent vectors \cite{eis}.    Conjugate nets parameterize constant mean curvature, constant Gauss curvature and isothermic surfaces, among others \cite{bob1}.   With this in mind, we consider the discretization of a parameterization on a surface by a ``discrete net'' of discrete coordinate curves: a map $f: \Z^2 \to \RP^n.$  Discrete surfaces in the classical geometries are then obtained by quadratic constraints on the image of the discrete net  in $\RP^n$ \cite{dol3}, following Klein's projective model  \cite{yag}. Thus, discrete conformal differential geometry, in the three-dimensional case, considers discrete nets in $\RP^4$  constrained to lie in the light-cone of Lorentzian metric $\R^5.$  Discrete nets are determined by solutions of systems of partial difference equations. In particular, the difference equations defining discrete conjugate nets \cite{dol1} are determined by requiring that the image of each face in $\Z^2$ lies in a projective plane. Discrete conjugate nets form the fundamental example of a discrete integrable system defining a discrete surface as solutions of the difference equation:


\begin{equation*}
\Delta_i\Delta_j f = a_{ji} \Delta_i f + b_{ij} \Delta_j f + c_{ij}f,
\end{equation*}  where $i,j$ vary over the lattice coordinates of $\Z^2$ \cite{bob1}. Results showing covergence from discrete to smooth solutions of the Laplace equation indicate that  the discrete and smooth surfaces are closely related \cite{bob3}.               

Just as Euclidean geometry is a subgeometry of M\"obius geometry,  Riemannian geometry is generalized by conformal (M\"obius) differential geometry \cite{mob}. M\"obius geometry is a subgeometry of projective geometry so that, in conformal differential geometry, parameterizations are coordinate maps from a surface into projective space. Isothermal coordinates, parameterizations compatible with the conformal structure of a surface, are employed as a natural construction.  Curvature-line coordinates, described typically in terms of the Riemannian principal curvatures end up being conformally invariant.   Thus, in conformal differential geometry it is natural to study isothermic surfaces: surfaces parameterized by curvature line coordinates compatible with the conformal structure of the surface. 

As the notion of a tangent space embedded at each point of the surface in ambient $3$-space is not conformally invariant it is natural to consider the conformally invariant ``mean curvature sphere'' \cite{h} of a surface.  This is a map assigning to each point of the surface in space a tangent sphere with radius given by the inverse of the mean curvature. Hence, the mean curvature sphere is properly a map from the surface into the space of two-spheres: a surface made of spheres.   Lie geometry, introduced by Sophus Lie in his doctoral dissertation under Pl\"ucker (and in collaboration with Felix Klein \cite{yag},)  is the geometry of oriented hyperspheres in M\"obius geometry.  In Klein's projective model, Lie geometry is the geometry of the Lie quadric, a signature $(2,4)$ real projective quadric. The set of null projective lines in the Lie quadric are an elementary example of contact geometry, with each null line corresponding to the one-parameter family of spheres tangent with a given orientation at a fixed point in space.  As the elementary invariant sets of M\"obius geometry are spheres, there is a natural relationship between conformal geometry and Lie geometry. Indeed, M\"obius geometry is contained as a subgeometry of Lie geometry with the $3$-sphere contained in the Lie quadric as the set of $2$-spheres of zero radius.  The M\"obius group is the subgroup of symmetries preserving these point-spheres, ``point transformations'' in the general  group of symmetries of the Lie quadric \cite{ww}.  A surface in space naturally induces a map into the space of contact lines in the Lie quadric determined by the set of tangent spheres at each point of the surface.  Each contact line contains a point corresponding to a sphere of $0$-radius, a point in $S^3.$ Thus, a map into the space of contact lines which satisfies the Legendre condition determines a surface in $S^3$ \cite{cec}.  

In the survey article of Bobenko and Suris \cite{bob2} the theory of discrete differential geometry is extended to contact lines in Lie geometry. Discrete ``principal contact element nets'' model the family of contact lines associated to a surface parameterized in curvature line coordinates.  The discrete analog of curvature line coordinates for surfaces in space are ``circular nets'', maps $\Z^2 \to S^3$ where the points of each elementary quadrilateral of the lattice lie on a circle.  Then, each contact line contains a representative corresponding to a point in $S^3$ so that the set of such point-spheres determines a circular net in $S^3$ \cite{bob1}. In the light-cone model of $S^3,$ a circle is given by the intersection of a projective $2$-plane with $S^3 \subset \RP^4.$  Thus, circular nets in $S^3$ are conjugate nets in $\RP^4$ subject to the constraint determined by the quadratic form of the light-cone \cite{dol3}.  

Circular nets in $S^2$ may be considered as a discretization of the Gauss map for a discrete surface \cite{bob3} \cite{udo2}. Identifying $S^2$ with $\CP^1,$ the cross-ratio of four circular points $[q_1,q_2,q_3,q_4] \in \R.$  Thus, circular nets in $S^2$ are solutions of  the ``cross ratio system,'' a set of partial difference equations given by determining the cross-ratio for each face of the map on $\Z^2.$    Given three points $\{q_1,q_2,q_3\} \in \CP^1,$ a complex number $\lambda$ determines a unique fourth point $q_4$ so that the four points have cross-ratio given by $\lambda.$  This cross-ratio system extends naturally to complex values and general discrete nets in $S^2.$ Thus in $S^2,$ the complex cross-ratio defines a master system of which circular nets are special solutions.

Let $c$ be  a discrete curve given as a map $c: \Z \to \C \subset \CP^1,$ then a choice of an initial point $c^{+}(0)$ defines the first iteration of the discrete evolution of $c$ by the formula:
\begin{equation*}
[c(1), c(0), c^{+}(0), c^{+}(1)] = \lambda.
\end{equation*}
Given $\lambda \in \C$ as a parameter, $c^{+}(1) = M_0(\lambda) c^{+}(0),$ where $M_0(\lambda)$ is a M\"obius transformation of $S^2.$ If $\lambda \in \R,$ then the points $\{c(1), c(0), c^{+}(0), c^{+}(1)\}$ are concircular.  Iterating this procedure, specifying $\lambda : \Z^2 \to \R$ determines a circular net in $\S^2.$ Thus, a discrete net may be viewed as the discrete evolution of a discrete coordinate curve. If $c$ is a closed curve, that is $c: \Z \to \H$ is periodic, then iterating $M_k(\lambda)$ around $c$ one is led to the holonomy problem  
\begin{equation}\label{curveflow}
 M_n(\lambda) ...M_1(\lambda) M_0(\lambda)c^{+}(0) =c^{+}(0)
\end{equation}   
given by the closing condition around $c.$  Thus, the eigenlines of the holonomy problem determine initial conditions for the evolution of a closed curve $c.$ 

Following interest in the physics of integrable systems such as the Novikov-Veselov equations, recent work has explored the conformal evolution of surfaces in four-dimensional M\"obius geometry \cite{sch}.  Concurrent work has discretized isothermic surfaces in $S^4$ by circular nets  \cite{pink}, \cite{udo2}.  We consider the quaternionic projective line $\HP^1 \cong S^4$ and then the quaternions  as coordinates of points in an affine chart  $ \H \subset \HP^1.$ We then identify $\R^3$  with the imaginary quaternions $ Im(\H) \subset \H \subset \HP^1.$   Given four points $q_1,q_2,q_3,q_4 \in \H$, by analogy with $\CP^1,$ it is possible to define a quaternionic cross-ratio by naively extending the complex formula to take quaternionic values relative to  an affine chart on $\HP^1$  
\begin{equation*}
[q_1,q_2,q_3,q_4] = (q_1 - q_2)(q_2 - q_3)^{-1}(q_3 - q_4)(q_4- q_1)^{-1} \in \H = \C \oplus j\C.
\end{equation*}
It is then possible to describe circular nets in $\HP^1$ as maps $f : \Z^2 \to \H$ where the quaternionic cross-ratio of the points in the image of each elementary face in $\Z^2$ is a real number. Discrete isothermic surfaces can then be defined in $\R^4 \cong \H \subset \HP^1$ by requiring the cross-ratio of each face to be a constant equal to $-1$ \cite{bobpink}.   Analogous conditions for discrete nets in $\C$ determine discrete holomorphic nets \cite{bob1} discretizing holomorphic functions.  

 Extending the algebraic integrable systems approach to smooth surface theory to quaternionic line bundles has provided powerful results related to geometric problems \cite{plu}.  It has been possible to extend the algebraic integrable systems approach to discrete problems.   In particular,  the cross-ratio evolution of a discrete closed curve in $\CP^1$  by (\ref{curveflow}) determines a discrete integrable system with ``spectral curve'' determined by the variety associated to the characteristic polynomial of the holonomy operator. The spectral parameter is given by the cross-ratio $\lambda,$ were $\lambda$ may take values in $\C$ \cite{per}.   For real values,  the quaternionic cross-ratio may be used to describe the discrete integrable evolution of discrete curves \cite{dol7} \cite{pink} in $S^4$.  Unfortunately,  the quaternionic cross-ratio is not a four-dimensional M\"obius invariant except for real values, i.e. for circular nets. Given four points in $\HP^1$ with quaternionic cross-ratio $q,$  it is possible to define a M\"obius invariant \emph{complex} cross-ratio by the formula $Re(q) + i|\Im(q)|$ \cite{bobpink}. However this complex cross-ratio does not determine a cross-ratio system: three points and a complex cross-ratio do not define a unique fourth point.   Thus, given the real cross-ratio system in $\HP^1,$ it was not clear how to extend the discrete integrable system to non-real values of the putative spectral parameter.

We propose a master discrete system which generalizes circular nets (real cross-ratio system) in $\HP^1$ and has a complex spectral parameter.  This system is defined by exploiting the flat twistor construction of Penrose.  With respect to the twistor projection from $\CP^3 \to \HP^1,$  each projective line in $\CP^3$ corresponds either to a point in $\HP^1$ or to a round two-sphere.  Those projective lines transverse to the twistor fibration correspond exactly to round two-spheres in $S^4$ \cite{h} \cite{bohle}. Hence, the set of round two-spheres and points in $\HP^1$ is identified with the complex Grassmannian  $G(2,4),$ itself considered as the four-dimensional Pl\"ucker quadric.  Thus, $\HP^1$ is embedded in the Pl\"ucker quadric as the light-cone model of the four-sphere.  This embedding is identified as the real set associated to a real structure on the Pl\"ucker quadric. Now, three non-collinear points in the Pl\"ucker quadric define a nondegenerate planar conic section.  We show that the Steiner cross-ratio \cite{dol1} on that conic section naturally extends the quaternionic  cross-ratio when those initial points lie in $S^4.$ Thus, the Steiner cross-ratio system on the Pl\"ucker quadric generalizes the cross-ratio system in $S^4$ where solutions typically consist of discrete nets of two-spheres in $S^4.$  In the smooth theory an example of such a surface made of spheres would be the map assigning to each point of an embedded surface the mean curvature sphere at that point.


 In this paper we have written an introduction to the elementary geometry of the complex Pl\"ucker quadric as ``twistor geometry:'' the geometry of two-spheres in the four-sphere.  This is a comprehensive and geometrically appealing interpretation to twistor theory that has, to our knowledge, not been explained elsewhere.  The geometry of the Pl\"ucker quadric thus stands as a neo-classical Kleinian geometry.  Indeed, we  show that the Lie quadric is contained as  a real subset of the Pl\"ucker quadric  with respect to a naturally defined real structure. The $(2,2)$ hermitian form on $\C^4$ of twistor theory is determined precisely by the choice of a three-sphere in $S^4.$ Thus, twistor geometry is obtained from the geometry of two-spheres in the three-sphere by complexification. In the literature of twistor theory, the four-sphere and the Lie quadric are identified as the conformal compactifications of Euclidean and Minkowski space respectively.  By this approach we can examine the geometry of two-spheres and circles in $S^4$ in terms of either projective lines (twistors) and (ruled) quadrics in $\CP^3$ or points and conic sections of the Pl\"ucker quadric, as convenient.  

In passing from codimension-$1$ to codimension-$2$ one obtains new contact properties for two-spheres.  We consider generalized discrete principal contact elements nets in the Pl\"ucker quadric and explain the novel geometry involved.  Complex null lines in the Pl\"ucker quadric consist of either a contact element of two-spheres tangent at a fixed point in $S^4$ or  a generalized contact element of two-spheres in ``half-contact'' at two points in $S^4.$  This idea was first introduced in \cite{pp}, determining the construction of generalized Darboux transforms of the torus in $S^4.$  Employing this complex contact geometry, we then show that the Steiner cross-ratio system contains both the real and complex cross-ratio systems in $\CP^1$ as  special solutions.  This leads to the novel result that the complex cross-ratio system in $\CP^1$ is in fact determined by a conjugate net system in a complex quadric.  Thus, discrete differential geometry in twistor space generalizes the theory of Bobenko and Suris for the Lie quadric.       


As the twistor viewpoint relies on easy switching between the natural objects of their respective spaces via equivalences, certain constructions in this paper are illustrated by a triptych of equivalent diagrams in the $4$-sphere, $\CP^3,$ and the Pl\"ucker quadric.  In particular, we indicate complex projective lines by pictures of configurations of lines in $\R^3$ and diagrammatically express configurations of points and $2$-spheres in $S^4.$  
          
\section{Twistor geometry: circles and two-spheres in $\HP^1.$}

The quaternionic model for conformal differential geometry in $S^4$ has been developed in \cite{h},\cite{mob}. Using the quaternions $\H$ instead of $\C$, the relationship between  $S^4$ and $\HP^1$ is entirely analogous to that between $\S^2$ and $\CP^1.$ Let $v \in \H^2,$ considered to be a two-dimensional quaternionic (right) vector space, then $v \mapsto [v] = v\H \in \HP^1$ defines the projective fibration of $\H^2$ over $\HP^1.$  Given a basis $\{e_1,e_2\}$ for $\H^2,$ $\HP^1$ is covered by affine charts $[e_1\H + e_2]$ and $[e_1 + e_2 \H]$ so that $\HP^1 = [e_1 \H + e_2] \cup [e_1] \cong \H \cup \{\infty\} = S^4.$ Orientation preserving M\"obius transformations are given by fractional linear maps induced by the natural left action of $\GL(2,\H)$ on $\H^2.$ We will consider, as M\"obius objects, points and round two-spheres.

Consider $\H$ as a real vector space spanned by $\{1,i,j,k\}$ with $i =jk$ and $i^2 = j^2 = k^2 = -1.$ Then, $\H = Re(\H) \oplus Im(\H)$ where $Im(\H)$ is given as the real span of the preferred orthonormal basis $\{i,j,k\}.$ Now, we may write $q \in \H$ in the form 
\begin{equation}
q = z_1 + j z_2,
\end{equation} 
for $z_1,z_2 \in \C.$  This may be seen by writing $a + ib + jc + kd = (a + ib) + j(c - id).$ By identifying the complex numbers with the real span of $\{1,i\} \subset \H,$ we then write $\H = \C \oplus j\C  \cong \C^2$ and $\H^2 \cong \C^4.$ If $\{e_1,e_2\}$ is a basis for $\H^2,$ then $\{e_1,e_1j,e_2,e_2j\}$ is a basis for $\C^4.$

\begin{Def}[Twistor Projection]
 Let $v \in \C^4 \cong \H^2$, then define the twistor projection $\tau : \CP^3 \to \HP^1$ by $v\C \mapsto v\H.$
\end{Def}

\begin{figure}[h]
\includegraphics{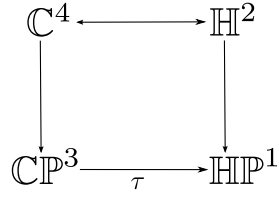}
\caption{The twistor diagram.}
\end{figure}

Let  $[v] \in \HP^1,$ then $[v] = v\H= v\C \oplus vj \C$ corresponds to a complex two-dimensional linear subspace of $\C^4$. Hence, the twistor fiber $\tau^{-1}([v])$  is a one dimensional projective line in $\CP^3$ spanned by $\{v, vj\} \subset \C^4.$ As each twistor fiber projects to a distinct point in $\HP^1$ each twistor fiber must be disjoint from every other twistor fiber. Thus, one may refer without ambiguity to the twistor fiber through a point $[v] \in \CP^3.$       

\begin{Def}
 A (round) two-sphere in $\HP^1$  is determined by the image under M\"obius transformation of any two-dimensional real subspace in an affine chart $\R^4 \cong \H \subset \HP^1.$
\end{Def}

In an affine chart on $\HP^1,$ a round two-sphere may be identified by M\"obius transformation with $\C \subset \H,$ with complex coordinates $z$ or $\overline z$ depending upon the conformal structure of the two-sphere. The following technical lemma will be used \cite{h}: 
\begin{Lem}\label{nr}\label{so3}
Let $v \in \H$ such that $|v| = 1,$ then $v \in \Im\H$ if and only if $v^2 = -1.$ 
\end{Lem}

Let  $\sigma \subset \HP^1$ be a two-sphere and $[\tilde e_1],[\tilde e_2]\in \sigma$ such that $\H^2 = \tilde e_1\H + \tilde e_2\H.$ Then, with respect to the affine chart $[\tilde e_1 \H + \tilde e_2],$ $\sigma$ is identified as a real two-dimensional subspace of $\R^4 \cong \H$ with $\infty = [\tilde e_1] \in \sigma,$ $0 = [\tilde e_2],$ and $1 = [\tilde e_1 + \tilde e_2].$  The two-sphere $\sigma$ may be identified with $\C\subset \H$ by scaling the original basis $\{\tilde e_1 \alpha,\tilde e_2 \beta\}$ for $\alpha, \beta \in \H.$ Suppose that $1 \notin  
\sigma.$ Given $[\tilde e_1 a + \tilde e_2]\in \sigma,$ then $\sigma$ contains $1$ in the affine chart $[\tilde e_1 a \H + \tilde e_2] \subset \HP^1.$  Extend $1$ to an orthonormal real basis $\{1,n\}$ for $\sigma.$  As $n$ is orthogonal to $1,$ $n \in Im(\H).$  Let $\lambda \in \H$ then $\lambda$ acts on $\R^3 \cong Im(\H)$ as an orientation preserving orthogonal transformation by $x \mapsto \lambda x \lambda^{-1}.$  Then, as $|n|=1$ and given Lemma \ref{so3} there exists $\lambda \in \H$ such that $n = \lambda i \lambda^{-1}.$  Thus, writing $e_1 = \tilde e_1 a \lambda$ and $e_2 = \tilde e_2 \lambda$ in the affine chart $[e_1 \H + e_2 ] \cong \H,$ $\sigma \setminus \{\infty\}$ is identified with $[e_1 \C + e_2] \cong \C.$

Define $(\sigma) \in End(\H^2)$ with respect to the basis $\{e_1,e_2\}$ by 
\begin{equation}\label{Sii}
 (\sigma) = \begin{pmatrix} i && 0 \\ 0 && i \end{pmatrix}.
\end{equation}
 Let $l = [e_1 x + e_2] \in \H \subset \HP^1$ such that $(\sigma)l = l.$ Then, 
\begin{equation}
 \begin{pmatrix} i && 0 \\ 0 && i \end{pmatrix} \begin{pmatrix} x \\ 1 \end{pmatrix} =  \begin{pmatrix} x \\ 1 \end{pmatrix} i
\end{equation}
which implies $xi = ix.$ That is, $x \in \C \subset \H$ and the set $\{l \in \HP^1$ such that $ (\sigma)l = l\}=\C\cup \{[e_1]\}.$   

Now, one may  express $(\sigma)$ in the original basis $\{\tilde e_1, \tilde e_2\}$ by
\begin{equation}\label{ana}
 (\sigma) = \begin{pmatrix}
             a n a^{-1} && 0 \\ 0 && n
            \end{pmatrix}.
\end{equation}
Hence, it is clear that the set of real two-planes containing $0$ in an affine chart $\H \subset \HP^1$ is determined exactly by elements of $\End(\H^2)$ of the form:
\begin{equation}
 \begin{pmatrix}
             R && 0 \\ 0 && N
            \end{pmatrix}
\end{equation}
for $R,N \in \H$ where $R^2 = N^2 = -1.$ Note that it is always possible to relate $R$ and $N$ by an orthogonal transformation of $Im(\H).$ 

Now, any two-sphere in $S^4$ may be obtained as the translation of a sphere through the origin. Let $[p]$ and $[f_2] \in \HP^1$ such that $p = e_1 q + e_2$ and $f_2 = e_2 + p.$ Suppose that $(\sigma ')e_1 = e_1 i$ and $(\sigma ')f_2 = f_2 i$ so that the two-sphere $\sigma '$ is given as the translation of $\C \cong [e_1\C + e_2]$ by $q.$ We compute 
\begin{equation}
(\sigma ') =  \begin{pmatrix}
             i &&  \frac{1}{2}(qi - iq) \\0 && i
            \end{pmatrix}.
\end{equation}  
Note that if $q \in \C$ then $\sigma = \sigma '.$  In general one obtains that each translation of $(\sigma)$ by  $q \in \H,$ is represented by 
\begin{equation}\label{RN}
\begin{pmatrix}
             R && H \\ 0 && N
            \end{pmatrix},
\end{equation} where $H = \frac{1}{2}(qN - Rq).$ If $q \in \sigma$ then the translation of $\sigma$ is trivial.  One obtains that $q \in \sigma \subset \H$ if and only if $qN - Rq =0$ and a proper translation $\sigma '$ of $\sigma$ is described by $[e_1 q + e_2] \in \sigma '$ if and only if $qN-Rq = 2H$ for some $H \in \H.$

If we consider the complex vector space $e_1 \C + e_2 \C \subset \C^4 \cong \H^2,$ then $\P(e_1 \C + e_2\C)$ is a projective line in $\CP^3.$ It is clear that $\tau(\P(e_1 \C + e_2\C)) = \sigma$  However, $\P(e_1 j \C + e_2 j\C) \cap \P(e_1 \C + e_2\C) = \emptyset$ and $\tau(\P(e_1 j \C + e_2 j\C)) = \sigma$  Then, as $End(\H^2) \subset End(\C^4),$ $(\sigma)$ acts on $\C^4.$ We may thus distinguish between the two previous sets by the property that $e_1\C + e_2\C \}$ is the $i$-eigenspace of $(\sigma).$  The unit quaternion $j$ acts on $\C^4$ by right multiplication. This action is anti-linear on $\C^4$ as $zj = j \bar z$ for $z \in \C.$ Thus, if $(\sigma)v=vi,$ 
then $(\sigma)vj=vj(-i)$ and $\{e_1 j\C + e_2 j\C \}$ is the $-i$-eigenspace of $(\sigma).$ We say that $\P(e_1 \C + e_2\C)$ is the twistor lift of $\sigma.$ When unambiguous, we will identify $\sigma$ with either its twistor lift or its image in an affine chart writing $\sigma = [e_1 \C + e_2].$   
Observe that 
\begin{equation}
\begin{split}[e_1 j z + e_2 j] = & [e_1 \bar{z} j + e_2 j] \\ = & [e_1 \bar{z} + e_2].
\end{split}
\end{equation}
Hence $\tau(\sigma)$ and $\tau(\sigma j)$ consist of the same set in $\HP^1$ with opposite conformal structures.  Note that two-spheres $\sigma$ and $\sigma '$ can be parallel with opposite conformal structures. If we write $\sigma =  [e_1 \C + e_2]$ and $\sigma ' =    [e_1 j \C + e_1 qj + e_2 j],$ then $\sigma'$ is a translation of $\sigma$ by $q \in \H$ followed by  conjugation.   
Thus, one has in \cite{h} the fundamental lemma:
\begin{Lem} \label{Stheory}
Let $S \in End(\H^2)$ such that  $S^2 = -Id$ then $S$ defines a unique round two-sphere in $\HP^1$ by $l\in \HP^1$ such that $Sl = l.$ For each two-sphere $\sigma \subset \HP^1$ with a given conformal structure, the twistor lift associates to $\sigma$ a unique projective line in $\CP^3.$
\end{Lem}
Now, let $l$ be a projective line in $\CP^3$ with the property that $l$ is disjoint from its image under the action of the quaternionic $j.$ Then $\{l, l j\}$ in $\CP^3$ corresponds via twistor projection to a unique two-sphere $\tau(l) \subset \HP^1.$  If $l = l j,$ then $l$ corresponds to a twistor fiber in $\CP^3$ that is, $\tau(l) \in \HP^1.$
Thus, one obtains: 
\begin{Pro}
Let $l \subset \CP^3$ be a projective line. Then, the image of $l$ under twistor projection is either a point in $\HP^1$ or a round two-sphere.
\end{Pro}

\begin{Rem}
Let $[x] = x\H \in \HP^1,$ then we can identify $Hom([x], \H^2 / [x]) \cong T_{[x]}\HP^1$ by $F \mapsto d_x\pi(F(x)$ where $\pi : \H^2 \to \HP^1.$
Given $\sigma,$ a two-sphere in $\HP^1,$ then $(\sigma)$ defines a complex structure on $T\sigma$ acting pointwise on $T_{[x]}\sigma \subset Hom([x], \H^2 / [x])$ \cite{h} . Hence, $(\sigma)$ defines a conformal structure on the two-sphere.  Given a surface $\Sigma$ embedded in $\HP^1,$ a two-sphere tangent  at each point of thus defines a conformal structure on $\Sigma.$    
\end{Rem}


\subsection{Twistor geometry}
 Let $Q^4$ be the complex quadric in $\P(\bigwedge^2 \C^4) \cong \CP^5$ defined by the zero set of the quadratic form induced by the wedge product $ \wedge: \bigwedge^2\C^4 \to  \bigwedge^4\C^4 \cong \C.$   Let $\alpha = v \wedge w \in \bigwedge^2\C^4 \cong \C^6,$ then associate to $[\alpha] \in \P(\bigwedge^2 \C^4)$ the projective line in $\CP^3$ spanned by $\{v,w\} \subset \C^4.$  Let $[\alpha] \in  \P(\bigwedge^2 \C^4),$ then $[\alpha]$ corresponds to a  line in $\CP^3$ if and only if $\alpha$ is decomposable, that is $\alpha = v \wedge w.$ By Poincare's Lemma, $\alpha$ is decomposable if and only if $\alpha \wedge \alpha = 0.$  Thus, there is an exact correspondence between $Q^4$ and  the set of projective lines in $\CP^3.$ We will use $[v \wedge w]$ to describe the point in $Q^4$ and its corresponding line in $\CP^3$ interchangeably. 

$Q^4$ is classically referred to as the Pl\"ucker quadric. In the twistor literature it is called the Klein quadric after Pl\"ucker's student, Felix Klein \cite{ww}. $ \C^4 \cong  \H^2 ,$ identified with a polar projection of  the Klein quadric, serves as the ``twistor space'' of Penrose's theory \cite{pen}. The Pl\"ucker quadric itself is interpreted as  ``complex compactified Minkowski space''  \cite{ww}.

\begin{Def}[Real Structure]
 Let $C$ be a complex manifold and $\rho: C \to C$ be a anti-holomorphic isomorphism. Then $\{\rho, C\}$ is a real structure.  The set of fixed points (possibly empty) of $\rho$ are said to be the ``real set`` of the real structure.  
\end{Def}

The right action of the scalar $j$ on $\H^2$ defines a right action on $\C^4 \cong \H^2.$  Following the twistor diagram, this defines a real structure on $\CP^3.$   The real set of $\{j, \CP^3\}$ is empty: no one-dimensional subspace of $\C^4$ is preserved by the right action of $j.$ However, this action may be defined on the set of projective lines in $\CP^3$ by extending anti-linearly the operation $(v \wedge w) j = vj \wedge wj$ on $\bigwedge^2 \C^4.$  It is well know in the literature of twistor theory \cite{ww} that this action of $j$ induces a real structure on $Q^4.$ The representatives in $Q^4$ of twistor fibers are given by points $[v\wedge vj] \in Q^4.$  These are exactly the fixed points of this real structure.  Hence, the real set of $j$ is equivalent to $\HP^1.$  In fact, it may be shown by direct computation that the real set of $j$ is given by a light-cone model of $S^4$ inside $Q^4,$ so that for $v \in \H^2$ the map $[v] \mapsto [v\wedge vj]$
determines the equivalence between $\HP^1$ and $S^4.$ 
In this fashion the geometry of the Pl\"ucker quadric may be said to be the complexification of four-dimensional M\"obius geometry. 
\begin{Pro}
Let $[\alpha] \in Q^4,$ then $[\alpha]$ is either a point in $\S^4 \subset Q^4$ or $[\alpha]$ is identified
with a two-sphere contained in $S^4$ and $[\alpha j]$ is identified with the same subset of $S^4$ with the opposite conformal structure. 
\end{Pro} 

We will now discuss the Pl\"ucker quadric in terms of the M\"obius geometry of two-spheres in $S^4.$

\begin{Def}[Touching Spheres]
 Let $\sigma$ and $\sigma '$ be two-spheres in $\HP^1$ each containing the point $[p] \in \HP^1.$ Then, $\sigma$ is \emph{tangent} to $\sigma '$ at $[p]$ if and only if, with respect to some basis $\{p,e_2\}$ spanning $\H^2,$ $\sigma$ and $\sigma '$ are given by parallel planes in the affine chart $[p\H + e_2] \subset \HP^1.$ We say that $\sigma$ \emph{touches} $\sigma '$ at $[p]$ if their respective conformal structures agree at $[p].$    
\end{Def}

The following is the fundamental relationship between the contact geometry of two-spheres in $S^4$ and the geometry of $Q^4.$

\begin{The}\label{RNT}
    Let $\sigma$ and $\sigma '$ be two-spheres  interesecting at a point $[p] \in \HP^1.$
Then, $\sigma$ and $\sigma '$ touch at $[p]$  if and only if the twistor lifts of $\sigma$ and $\sigma '$ are incident  and the twistor fiber over $[p]$ lies in the plane defined by $\sigma$ and $\sigma '.$ 
\end{The}

{\em Proof:} Let $\sigma$ and $\sigma '$ be two-spheres touching at a point $[p] \in \HP^1.$  Then, $\sigma$ and $\sigma '$ touch at $[p]$ if and only if their representatives in $End(\H^2),$ given in the form of (\ref{RN}), share the same $R$ and $N.$  By construction $(\sigma) p = pR = (\sigma ') p$.   Given $R^2 = -1,$ using Lemma \ref{so3}  there exists $\lambda \in \H$ such that $\lambda^{-1} R \lambda = i.$ Then,
\begin{equation}
 (\sigma) p\lambda = p R \lambda = p \lambda (\lambda^{-1} R \lambda) = p \lambda i.
\end{equation}
Hence the twistor lifts of $\sigma,$ $\sigma '$ and $[p]$ in $\HP^1$ intersect at the point $[p\lambda] \in \CP^3.$ 

Let $\{p,q\}$ be a basis for $\H^2.$ Consider the left quaternionic dual space ${\H^2}^*$ of linear quaternionic functions on $\H^2.$ We may now consider an analogous twistor construction by identifying ${\H^2}^*$ with ${\C^4}^*.$ Then $(\sigma)$ acts on ${\H^2}^*$ by precomposition leading to twistor lifts into ${\CP^3}^* = \P({\C^4}^*).$ Let $\pi \in p^o = \{ \alpha \in  {\H^2}^*$ such that $\alpha(p) = 0\}$ such that  $\pi(q)=1.$ Then,
\begin{equation}
 (\sigma)^* \pi = N \pi.
\end{equation}
Thus, we have shown that the dual twistor lifts of $\sigma$ and $\sigma '$ are incident at a point in the dual space to $\CP^3.$

\begin{figure}[h]
  \includegraphics{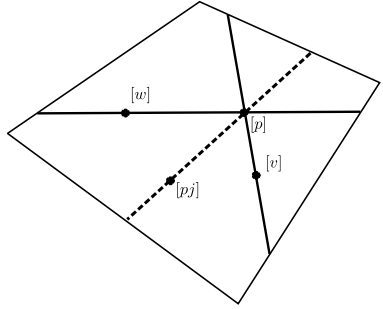}
\caption{Two-spheres touching at point $p,$ represented in incident lines in $\CP^3.$ The dotted line is the twistor fiber of $[p].$}
\end{figure}

By projective duality, if the twistor lifts of $\sigma,$ $\sigma '$ and $[p]$  intersect in a point in $\CP^3,$ then their dual representatives must be coplanar in the dual space. But then, we have just seen that these representatives also intersect in a point in the dual space. Reversing duality, we obtain that the twistor lifts of $\sigma,$ $\sigma '$ and the twistor fiber through $[p]$ must be coplanar in $\CP^3.$  

Finally, we may write $\sigma = [p\C + v]$ and $\sigma ' = [p\C + w].$ As the twistor fiber of $[p]$ is the projective span of $\{[p],[pj]\}$ and is contained the plane generated by $\{[p],[v],[w]\},$ we may write $w = p z_1 + pj z_2 + v,$ for $z_1,z_2 \in \C.$ Then, $w = pq + v$ for $q = z_1 + jz_2 \in \H$ and $\sigma ' = [p(\C + q) + v].$ Thus, $\sigma'$ is a translation  of $\sigma$ and their conformal structures agree.$\square$    

Given Theorem \ref{RNT}, we now recognize the representative of a two-sphere touching a given two-sphere $\sigma = [p\C + v]$ at the point $[p \wedge pj]$ as a line incident to $[p] \in \CP^3$ and coplanar with  $[p\C + v].$ The set of two-spheres touching $\sigma$ at $[p]$ is given by the complex one-parameter family of lines incident to the twistor fiber  $\tau^{-1}[p]$ at $[p]$ and contained in the projective plane spanned by $\sigma$ and $\tau^{-1}[p].$  We will now show that this family is paramterized by the points of a projective line contained in $Q^4.$
  
A linear subset of a projective quadric is said to be a null object, e.g. null line, null plane. Two linear subsets are said to be incident if they intersect. 
As $Q^4$ is a hypersurface in $\CP^5,$ the projective dimension of the maximal null object is two. We now state some standard lemmas 
describing the null objects in $Q^4$ \cite{ww}. As a line is defined by two points, the fundamental lemma is the following:
\begin{Lem}
 Let $[\alpha],[\beta] \in Q^4.$ Then $\alpha \wedge \beta = 0$ if and only if the projective lines $\alpha$ and $\beta$ are incident in $\CP^3.$ 
\end{Lem}

Thus, a null line in $Q^4$ is determined by two incident but not identical lines in $\CP^3.$ Two such incident lines in $\CP^3$ must lie in a projective plane and intersect in a point.

\begin{Lem}
A null plane in $Q^4$ corresponds to the set of projective lines in $\CP^3$ incident to a given point (type A) or, to the 
set of projective lines all contained in a given plane (type B).  
\end{Lem}

The set of lines contained in a plane in $\CP^3$ is, by projective duality, equivalent to the set of lines incident to a point.  Thus, the set of elements in $\CP^3$ corresponding to a null plane of type B in $Q^4$ is the image by projective duality of the elements of a null plane of type A in the dual space to $\CP^3.$

\begin{Lem}
Each null line is given by the intersection of null planes of type A and type B.
\end{Lem}

\begin{Cor} \label{pointplane}
Each null line corresponds to the pencil of lines in $\CP^3$ contained in a given plane and incident to a fixed point contained in that plane.  The set of null lines in $Q^4$ corresponds exactly to the set of pairs $\{\Pi, [p]\}$ where $\Pi$ is a projective plane in $\CP^3$ and $[p] \in \Pi.$ The set of points on that null line corresponds to the pencil of lines in $\CP^3$ incident to $[p]$ and contained in $\Pi.$
\end{Cor}

Now, the unit quaternion $j$ acts on $\C^4$ anti-linearly. Thus, it descends to a map on $\CP^3.$ In particular, the image of a projective plane in $\CP^3$ under $j$ is a projective plane. But, two projective planes in $\CP^3$ must intersect in a line. Hence,

\begin{Pro}\label{utf}
Each projective plane in $\CP^3$ contains exactly one twistor fiber. 
\end{Pro}
Two points are contained in a null line if and only if their representatives in $\CP^3$ intersect. Two lines intersect in $\CP^3$ if and only if they are mutually contained in a projective plane. Each null line in the Pl\"ucker quadric may be classified by whether it contains a point corresponding to a twistor fiber or not.  Recal that $\{Q^4, j\}$ is a real structure.

\begin{Cor}\label{classinull}A null line in $Q^4$ either contains exactly one real point of $j$ or no real points.
 \end{Cor}

Thus, we obtain, as a corollary to Theorem \ref{RNT}, that the set of two-spheres touching a given two-sphere $\sigma$ at a point $[p],$ including $\sigma$ and the point of tangency is parameterized by  a null line in $Q^4.$  

\begin{Pro}\label{touchingnull}
 Each null line in $Q^4$ containing a real point corresponds to the set of two-spheres in $\HP^1$ mutually touching at that point.
\end{Pro}

\begin{figure}[h]
 \includegraphics{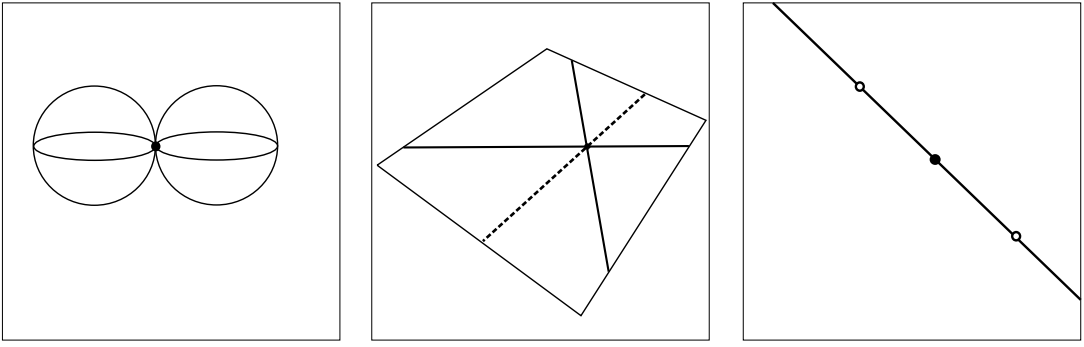}
\caption{Touching spheres and their point of contact in $\HP^1,$ $\CP^3,$ and $Q^4$}
\end{figure}

It is clear that not every null line in $Q^4$ contains a real point, that is, corresponds to a family of touching spheres.  A null line of spheres touching a given sphere $\sigma$ at $[p] \in S^4 \subset Q^4$ is defined by choosing a projective plane $ \Pi \subset \CP^3$ containing the twistor fiber $\tau^{-1}[p]$ and $\sigma.$ 
Now, each projective plane containing $\sigma$ must contain the twistor fiber $\tau^{-1}[q]$ of a point $[q]$ on $\sigma.$  Thus, we may deform $\Pi$ by the projective plane $\Pi_q$ defined by the projective span of $\sigma$ and $\tau^{-1}[q]$ for some $[q] \in \sigma,$ where $\Pi = \Pi_p.$

For each $\Pi_q$ we can define a null line in $Q^4$ by taking the set of projective lines incident to $[p] \in \Pi_q$ and contained in $\Pi_q.$  An element $\sigma '$ of this null line is thus given by $\sigma ' = [p \wedge q\mu]$ for some $q\mu \in \tau^{-1}[q].$  It is clear from Theorem \ref{RNT}  that this null line cannot consist of touching spheres as $\tau^{-1}[p]$ is not contained in $\Pi_q.$
Suppose that $(\sigma) p = pi$ and $(\sigma)q = qi.$  Looking at (\ref{RN}), we observe that $\sigma'$ has representative $(\sigma ') \in End(\H^2),$ with respect to the basis $\H^2 = p \H + q \H,$ of the form
\begin{equation}
 \begin{pmatrix} \mu i \mu^{-1} && 0 \\ 0 && i \end{pmatrix}.
\end{equation}
This is clear as $\sigma \cap \sigma ' = [p]$ and by construction $(\sigma ') q\mu = q\mu i.$

 \begin{figure}[h]
   \centering
 \includegraphics{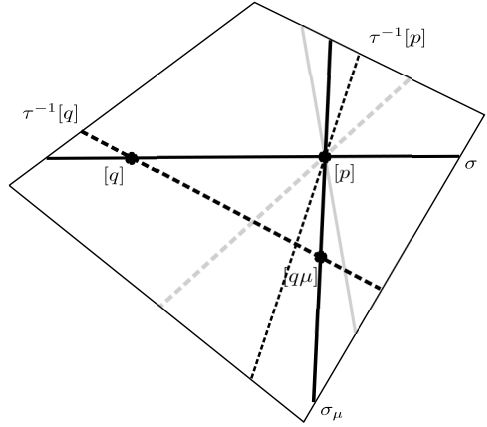}  
   \caption{Deforming $\Pi$ to $\Pi_q$ to obtain a null line of spheres which are not touching.}
 \end{figure}

Consider the set of two-spheres in $\HP^1$ containing the representatives of the real points $[p \wedge pj], [q \wedge qj] \in Q^4.$
The corresponding set in $Q^4$ is given by the nondegenerate complex two-dimensional quadric $Q_{p,q} = \{[p \wedge pj],[q \wedge qj]\}^{\perp} = \{ [a \wedge b] \in Q^4$  such that $ a \wedge b \wedge p \wedge pj  = a \wedge b  \wedge q \wedge qj = 0\} \subset Q^4.$ Let $\sigma \in Q_{p,q},$ then $\sigma^{\perp} \cap  Q_{p,q}$ consists of two null lines intersecting exactly at $\sigma.$  By construction neither of these null lines can contain a real point else that point would have, as a representative in $\CP^3,$ a twistor fiber intersecting both $\tau^{-1}[p]$ and $\tau^{-1}[q].$ If we examine this configuration in $\CP^3$ then $\tau^{-1}[p]$ and $\tau^{-1}[q]$ are disjoint projective lines with the line $\sigma$ incident to both. Suppose $\sigma \cap  \tau^{-1}[p] = [p\mu] \in \CP^3$ and $\sigma \cap  \tau^{-1}[q] = [q\nu]$ for $\mu, \nu \in \H.$ Write $(\sigma) \in End(\H^2)$ in the basis $\H^2 = p \H + q \H,$ as
\begin{equation} \label{rn}
 \begin{pmatrix} R && 0 \\ 0 && N \end{pmatrix}.
\end{equation}
Then, the set of lines incident to $[p\mu]$ and containing $\tau^{-1}[q]$ and the set of lines in $\CP^3$ incident to $[q\nu]$ and containing $\tau^{-1}[p]$ lie in unique projective planes $\Pi_\mu$ and $\Pi_\nu$ respectively. It is clear  that these families of lines define the two null lines incident at $[\sigma]$ and contained in $ \{[p],[q]\}^{\perp}.$ The spheres in the null line defined by $\{\Pi_\mu, [q\nu]\}$ are expressed in the form
  \begin{equation}
 \begin{pmatrix} R && 0 \\ 0 && * \end{pmatrix}
\end{equation}
and those in the null line defined by $\{\Pi_\nu, [p\mu]\}$ are expressed in the form
   \begin{equation}
 \begin{pmatrix} * && 0 \\ 0 && N \end{pmatrix}.
\end{equation}
Further, if any twistor fiber is incident to two members of either family then it must lie in the respective plane and thus be either $\tau^{-1}[p]$ or $\tau^{-1}[q].$

Thus, we are led to a concept introduced by Pedit and Pinkall in \cite{pp} that is central to the geometry of  two-spheres in four dimensional M\"obius geometry.

\begin{figure}[h]
  \centering
  \includegraphics{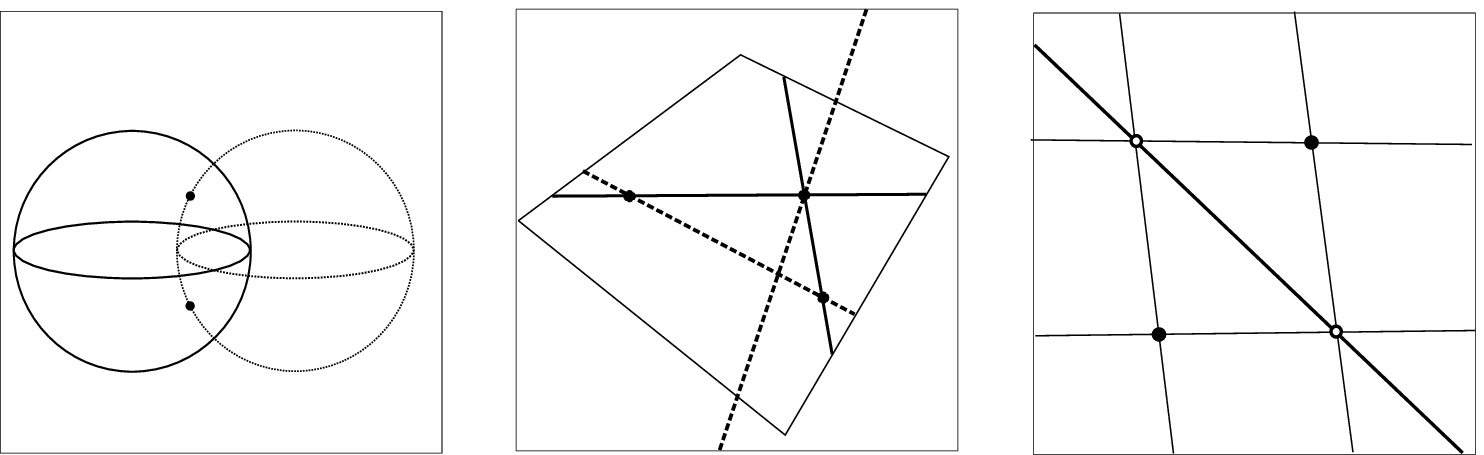}
  \caption{Half-touching spheres in $\HP^1,$ $\CP^3,$ and $Q^4.$}
\end{figure}

\begin{Def}[Half-Touching Spheres]
The family of spheres contained in a null line in $Q^4$ with no real point with respect to the real structure $j$ are said to be ``half-touching'' and intersect at exactly two points in $S^4.$ A family of spheres half-touching at $[p], [q] \in \HP^1$ are ``right-touching'' at $[p]$ if their twistor lifts in $\CP^3$ are all incident at some point on $\tau^{-1}[p].$ If they are right-touching at $[p]$ then they are said to be ``left-touching'' at $[q].$   
\end{Def}


 The tangent space to a smooth surface embedded in $S^4 \cong \HP^1$ is not M\"obius invariant.  However, given a two-sphere tangent to this surface at a point, writing the two-sphere in the form (\ref{rn}), the pair $\{R,N\}$ defines the conformal structure of the surface at that point.  Thus, two half-touching two-spheres share half of their respective conformal structures.  Two surfaces which intersect at an isolated point and with half-touching tangent spheres can said to be in ``half-contact.''  This extends the contact geometry of surfaces in $S^3$ to the codimension-$2$ contact geometry in $S^4.$   Contact elements in this geometry are null lines in $Q^4$ corresponding to either touching or half-touching spheres in $S^4$ according to the classification of null-lines given in Cor. \ref{classinull}. 


Now, let $\Pi = \{\sigma, \tau^{-1}[p]\}$ be the plane in $\CP^3$ defined by a twistor fiber $\tau^{-1}[p]$ and a two-sphere $\sigma = [p \wedge q]$ incident at $[p].$  We may now deform $\sigma \mapsto \sigma_{q\mu}$ within the family of half-touching spheres at $[p]$ by $[p \wedge q] \mapsto [p \wedge q \mu]$ for $\mu \in \H.$ By the previous discussion, $\sigma_{\mu}$ right-touches $\sigma$ at $[p].$ 

Let $\sigma = [p\C + q] \subset \CP^3.$ Writing $\H^2 = p\H + q\H,$ $\sigma \cong \C \subset \H.$ Then, $\sigma_{\mu} \subset \CP^3$ can be written
\begin{equation}
\sigma_{\mu} = [p \C + q\mu] = [p \C \mu^{-1} + q].
\end{equation}
In a similar fashion, we can deform $\sigma$ along $\tau^{-1}[p]$ by
\begin{equation}
 \sigma_{\lambda} = [p(\lambda  \C) + q] 
\end{equation}
so that $\sigma_\lambda$ left-touches $\sigma$ at $[p].$  Thus, $\sigma_\mu \cong \C \mu^{-1} \subset \H,$ and $\sigma_\lambda \cong \lambda \C \subset \H.$ It is clear then that given $[p],[q] \in \sigma$ an arbitrary two-sphere containing $[p]$ and $[q]$ is defined by
\begin{equation}
\sigma_{(\lambda,\mu)}= [p\lambda \C \mu^{-1} + q].
\end{equation}
This may be illustrated with the help of the following lemma \cite{bohle}:

\begin{Lem}\label{so4}
Let $\lambda, \mu \in \H$ such that $|\lambda| = |\mu| = 1,$ then $(\lambda, \mu) \in SU_2 \times SU_2,$ acts on $x \in \H$ as an orientation preserving orthogonal transformation of $\R^4$  by  $\lambda x \mu^{-1}.$  
\end{Lem}

For the sake of illustration consider 
\begin{equation}
 S = [p(j\C) + q],
\end{equation}
so that $S = j\C \subset Im(\H) \cong \R^3.$ Now, by Lemma \ref{so3}, if $\lambda \in \H$ such that $|\lambda|=1,$  then $\lambda j\C \lambda^{-1} \subset Im(\H)$ is a real plane containing $0.$ Then, by Lemma \ref{so4} we see that the diagonal group of $SU_2 \times SU_2$ acts by $SO_3$ on $\Im(\H).$  Thus the orbit of $j\C$ under $j\C \mapsto \lambda j\C \lambda^{-1}$ consists of all $2$-dimensional subspaces of $\R^3.$ Acting by $SO_4$ on $S$ using the left and right components of $SU_2 \times SU_2$  obtains the set of two-spheres in $S^4$ half-touching $S$ at $0$ and $\infty.$

\subsection{Circles}

A circle in M\"obius geometry is defined uniquely by three points. You can see this in $\HP^1$ in the fact that, for three points $\{[p],[q],[r]\},$ each affine chart with $[p] = \infty,$ $[q] = 0,$ and $[r] = 1$ identifies $\R \subset \H$ with $[p\R +q].$  We observe that three points in $\HP^1$ correspond to three skew projective lines $\{[p\wedge pj], [q \wedge qj], [r \wedge rj]\}$ in $\CP^3.$    A projective line incident to      
each of the twistor fibers $\{[p\wedge pj], [q \wedge qj], [r \wedge rj]\}$ corresponds to a two-sphere containing the three points and hence containing the associated circle in $\HP^1.$ Thus, the set of lines incident to the three twistor fibers corresponds to the set of two-spheres in $\HP^1$ containing the associated circle.

Using classical terminology, three skew lines in space define a ``regulus:'' a unique one-parameter family of mutually skew lines containing the initial three lines.  This family of lines sweeps out a complex quadric in $\CP^3$ \cite{ped}. The initial three skew lines are \emph{generators} of the regulus. The quadric may be constructed by observing that a point and disjoint projective line in space span a unique projective plane. Then, a line not contained in that plane must intersect it in exactly one point. Thus,  for each point on one of three skew lines there is a unique line incident to it and the other two lines.  A nondegenerate quadric in $\CP^3$ is ruled by two one-parameter families of mutually skew lines such that each line in one family intersects every line in the other family. Thus, the regulus is defined as the family of lines containing the original three skew lines and we say the other family  is ``polar'' to it . We say a line incident to three skew lines is incident to the regulus of those lines. We say that a one-dimensional projective (sub)quadric is a \emph{conic section} or conic.  

The regulus and its  polar may be represented as disjoint conics in $Q^4$ which are polar in the sense that each point in one conic is orthogonal to every point in the other.  The regulus of $\{[p\wedge pj], [q \wedge qj], [r \wedge rj]\}$ is defined by the intersection of the projective plane spanned by $\{[p\wedge pj], [q \wedge qj], [r \wedge rj]\}$ with $Q^4.$ This defines a conic section, each point of which corresponds to a line of the regulus in $\CP^3.$ The set of lines in $\CP^3$ incident to $\{[p\wedge pj], [q \wedge qj], [r \wedge rj]\}$ is defined by the intersection of the projective plane $\{[p\wedge pj], [q \wedge qj], [r \wedge rj]\}^{\perp}$ with $Q^4.$ This defines a second conic section of $Q^4$ and thus determines a one-parameter family of lines in $\CP^3.$  The real points of the regulus will correspond to the common points of intersection of a set of two-spheres each mutually containing three fixed points, that is, a circle.  We will compute explicitly the coordinates of points on the regulus and come to a geometric interpretation of the non-real points.     

\begin{figure}[h]
 \includegraphics{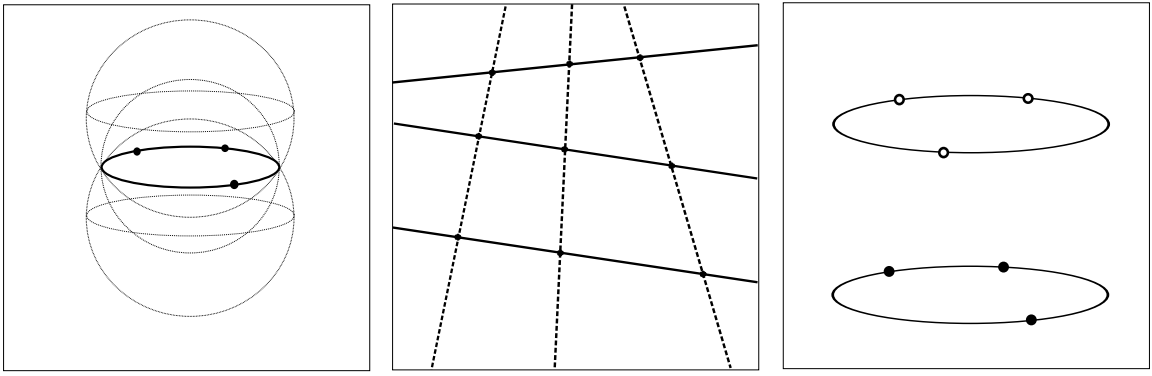}
\caption{The family of of two-spheres incident to a circle in $\HP^1$, the family as regulus of the circle in $\CP^3$, and the conic of the circle and regulus in $Q^4.$  }
\end{figure} 

\begin{Lem}\label{l1}
Let $\pi, \xi,$ and $\rho$ be three skew lines in $\CP^3$ and $S$ and $\tilde S$ lines incident to each of them. Consider $\alpha \in C_r \subset Q^4,$ the regulus associated to $\{\pi, \xi, \rho\},$  then
\begin{equation}
  \label{eq:regcoor}
  \alpha = [(p z + q)\wedge (\tilde p z + \tilde q)], 
\end{equation}
where $z \in \C,$ $p,q \in S,$ and $\tilde p, \tilde q \in \tilde S.$ 
\end{Lem}
{\em Proof}: Given $[\pi], [\xi],$ and $[\rho]$ there exist skew lines $S$ and $\tilde S$ incident to each of them.  Then 
$[\pi] = [p \wedge \tilde p], [\xi] = [q \wedge \tilde q],$ and $[\rho] = [r \wedge \tilde r]$ where $[p] \in S,
[\tilde p] \in \tilde S,$ and likewise for the others. Note that $\{ p,\tilde p, q, \tilde q \}$ span $\C^4$. Without 
loss of generality we assume that $r = p + q$ and $\tilde r = \tilde p + \tilde q.$ Now, let $\alpha \in \Span \{ 
\pi, \xi,\rho \}$ so that $\alpha \wedge \alpha = 0.$ Write
\begin{equation}
\alpha  =  r \wedge \tilde r + p \wedge \tilde p a+  q \wedge \tilde q b
\end{equation}
for $a,b \in \C.$ Then,
\begin{equation}
\begin{split}
  \alpha & = (p (1+a) + q)\wedge \tilde p + (q(1+b) + p)\wedge \tilde q \\
       & = \beta + w \wedge \tilde q,
\end{split}
\end{equation}
where $w = q(1+b) + p \in W = \Span\{p,q,\tilde p\}$ and $\beta \in \bigwedge^2 W$ so that $\alpha \wedge \alpha = 2\beta \wedge w \wedge \tilde q.$ But, $\tilde q \notin W$ and $\alpha \wedge \alpha = 0$ imply that $\beta \wedge w = 0 = (p(1+a) + q)\wedge w.$ Hence, 
\begin{equation}
 w = (p(1+a)  + q)k.
\end{equation}
Then $k(1+a) = 1$ and $k = 1+b$ so that
\begin{equation}\label{conic}
\begin{split}
 \alpha & = w \wedge ( \tilde p k^{-1}+ \tilde q) \\
 & = (p(1+a)  + q) \wedge ( \tilde p (1 + a)+ \tilde q). \square
\end{split}
\end{equation} 
We have now an explicit formula for the points of the regulus of three skew lines excluding the point $\pi.$ Each point on this curve corresponds to a projective line skew to $[\pi],[\xi],$ and $[\rho]$ which are themselves elements of the regulus.

Now, let $\tilde p = pj$ and $\tilde q = qj$ so that $\pi, \xi,$ and $\rho$ are three real points in $Q^4.$ By Lemma \ref{l1},$\pi, \xi,$ and $\rho$ are elements of a one-parameter complex family of lines in $\CP^3.$  Let $\sigma \in  C^{\perp}_r = \{\pi, \xi,\rho\}^{\perp} \cap Q^4,$ then $\sigma j \in C$ so that $\{\sigma,\sigma j\}$ corresponds to a two-sphere containing this circle. Let $\alpha \in C_r,$ then $\alpha$ corresponds to a line in $\CP^3$ incident to the projective lines corresponding to $\sigma$ and $\sigma j.$  If $\alpha$ is not a twistor fiber, then it is clear that the two-sphere represented by $\alpha$ must half-touch that of $\sigma.$ 

\begin{Cor}
Let $S= [p\C + q] \cong \C \subset \H$ so that the circle defined by the points $[p],[q],$ and $[r]=[p + q]$ in $\HP^1$ corresponds to $[p\R + q] \cong \R \subset \H.$ The twistor fibers $\{[p\wedge pj], [q \wedge qj], [r \wedge rj]\}$ define a conic section $C_r \subset Q^4.$ Points on $C_r$ correspond to either points on the circle $[p\R + q]$ including $\infty$ or two-spheres in $\HP^1$ which half-touch $S$ with respect to the given coordinates at pairs of points $\{z, \bar z\}$ in $\C.$   
\end{Cor}

{\em Proof}: Let $\alpha$ be given as in Lemma \ref{l1} for $\tilde p = pj, \tilde q = qj,$ and $\tilde S = Sj$:
\begin{equation}
 \alpha = [(pz  + q) \wedge ( pj z+ qj)].
\end{equation}
Then, the line in $\CP^3$ corresponding to $\alpha$ must intersect $S$ at $[pz + q]$ and $Sj$ at $[pjz + qj]$ for $z \in \C.$ But the point in $\HP^1$ corresponds to the twistor fiber of  $[pjz + qj]$ is given by
\begin{equation}
 [pjz + qj] = [(p\bar z + q)] = [p \bar z + q].
\end{equation}
since $jz = \bar z j.$ Thus, the two-sphere in $\HP^1$ represented by $[\alpha]$ must half-touch $S \cong \C \cup \{\infty\}$ at $z$ and $\bar z.$ $\square$

It is possible to project any conic $C$ onto a fixed line $l$ coplanar to it by fixing one point on the conic $p \notin l.$ Let $x \in C,$ then $x \mapsto \phi_l(x) \in l$ where x is mapped to the intersection of the projective line spanned by $\{x,p\}$ and $l.$  The point $p$ is then mapped to the intersection of its polar line and $l.$  Consider the regulus defined by $\{[p\wedge \tilde p],[q\wedge \tilde q],[r\wedge \tilde r]\},$ where $p,q,r \in S$ and $\tilde p,\tilde q ,\tilde r \in \tilde S$ for $S,\tilde S$ incident to the regulus. The points of $S$ parameterize the regulus by identifying each line in the regulus with it's intersection with $S.$  We will now  compute this parameterization explicitly for the regulus  $C_r$ and show that it is compatible with the map $\phi_l.$ 

\begin{Lem}\label{l2}
 Let $\{[p\wedge \tilde p],[q\wedge \tilde q],[r\wedge \tilde r]\},$ be generators of the regulus $C_r \subset Q^4$ and $S,\tilde S\ \in C^\perp.$ Given $[\eta \wedge \tilde \eta] \in C_r,$ then there is $x \in \C$ such that $[\eta] = [px + q] \in S$ and $[\tilde \eta] = [\tilde p x + \tilde q] \in \tilde S.$ Let $l_s$ be the secant line to $C_r$ defined by the $\Span\{[p\wedge \tilde p],[q\wedge \tilde q]\},$ then there exist coordinates on $l_s$ so that $x \mapsto \phi_l([\eta\wedge \tilde \eta])$ is the identity map.   
\end{Lem}
{\em Proof}: As in Lemma \ref{l1}, we may consider without loss of generality, 
\begin{equation}
\begin{split}
r\wedge \tilde r & = (p + q)\wedge(\tilde p + \tilde q)\\
& = p \wedge \tilde p + (p \wedge \tilde q +  q \wedge \tilde p) + q \wedge \tilde q.
\end{split}
\end{equation}
where $p,q \in S$ and $\tilde p, \tilde q \in \tilde S.$ 
Then, by Lemma \ref{l1}, $[\eta \wedge \tilde \eta] \in C_r$ is given by 
\begin{equation}
(p x + q) \wedge \tilde (p x + \tilde q) =   p \wedge \tilde p x^2 + (p \wedge \tilde q +  q \wedge \tilde p)x + q \wedge \tilde q.
\end{equation}
Now, fix $[r\wedge \tilde r] \in C_r$ and let $l_s$ be the secant line defined by the $\{[p\wedge \tilde p],[q\wedge \tilde q]\}.$ The intersection of the line defined by the $\{[r\wedge \tilde r], [\eta \wedge \tilde \eta]\}$ is thus given by
\begin{equation}
\begin{split}
 \phi_{l_s}([\eta\wedge \tilde \eta]) & = [\eta\wedge \tilde \eta - [r\wedge \tilde r]x]\\
 & = [p \wedge \tilde p x + q \wedge \tilde q]. \square
\end{split}
\end{equation}

Considering the map $\phi_l(x)$ as defining coordinates on a conic $C,$ then choosing another line $l'$ transforms these coordinates by a projective map. These coordinates are thus projective coordinates. Then, up to projective equivalence, the map defined by Lemma \ref{l1} from $S \to C$ is the identity map. Thus, Lemma \ref{l1} indicates that any line incident to a regulus defines projective coordinates for the conic in $Q^4$ associated to the regulus.       

\begin{Def}
 The Steiner cross-ratio of four points on a nondegenerate conic $C$ is a projective invariant of $C$ defined by the cross-ratio of those four points with respect to their projection onto some projective line coplanar to $C$. 
\end{Def}

Thus we have the following corollary to Lemma \ref{l2}:
\begin{Cor}\label{steiner}
 Let $S \subset \CP^3$ be any projective line incident to the regulus $C_r \subset Q^4$ of the skew lines $\{\pi,\xi,\rho\}.$    The Steiner cross-ratio of four-points on $C_r$ may be computed by taking the cross-ratio with respect to their coordinates on $S.$ 
\end{Cor}
Given the previous lemmas we now show that the Steiner cross-ratio on $Q^4$ generalizes the quaternionic cross-ratio.

\subsection{The Generalized Cross-Ratio}

It is well known that a choice of three points in $\CP^1$ defines unique projective coordinates up to an ordering of the points. With respect to these coordinates the initial three points may be chosen to correspond to $\{\infty,0,1\}$ and any other point $p \in \CP^1 \setminus\{\infty,0,1 \}$  is given in homogenous coordinates of the form 
\begin{equation}
p = \begin{bmatrix} z \\ 1\end{bmatrix}
\end{equation} with coordinate $z \in \C.$  
\begin{Def} Let $\{p_1,p_2,p_3,p_4\}$ be distinct, ordered points in $\CP^1.$  In the unique coordinates such that $\{p_1,p_2,p_3\} \to  \{\infty,1,0\}$  the  cross-ratio $[p_1,p_2,p_3,p_4]$ is defined as  the coordinate of $p_4$.  In arbitrary coordinates, we will identify points with their coordinates and write the cross-ratio of the four points  as  \begin{equation}[z_1, z_2, z_3,z_4 ] = \frac{(z_1 - z_2)}{(z_2 - z_3)}\frac{(z_3 -z_4)}{(z_4 -z_1)}. \end{equation} 
\end{Def}
Thus, in particular, $[\infty, 1,0,\lambda] = \lambda.$ Allowing the cross-ratio to vary freely over the real numbers, the points $\lambda$ determine the real line $\R \subset \C.$  The cross-ratio is invariant under $\GL_2(\C)$ acting by M\"obius transformations on $\CP^1$ and hence a conformal invariant.  Any permutation of $\{\infty, 1, 0, \lambda\}$ results in a number  which is a fractional linear transformation of $\lambda.$  Recognizing that any circle in $\C$ may be identified with the real line by M\"obius transformation, four points with real cross-ratio must lie on a unique circle or affine line  contained in $\C.$       

Thus, the classical cross-ratio of four points on the complex projective line $\CP^1$ has two elementary properties: 
\begin{itemize}
\item up to a permutation of order of the points, it is invariant  under conformal maps  on $S^2 \cong \CP^1$ 
\item three points and a complex number uniquely define a fourth point so that, with  respect to an ordering of the three initial points, the fourth point is given as the image by a unique M\"obius transformation of the third point and has cross-ratio given by the chosen number.
\end{itemize}

Defining $\HP^1$ by analogy to $\CP^1,$ replacing the field of complex numbers $\C$ with the the skew-field of quaternions $\H,$ it is natural to attempt to extend the cross-ratio to $\H$ by the previous definition.    Hence, given  three points $\alpha, \beta, \gamma  \in \HP^1$ and choose $c \in \H^2$ so that $\gamma = c\H.$ Now, there exist $a, b \in \H^2$ so that $\alpha = a\H, \beta = b\H,$ and $c = a + b \in \H^2 = a\H \oplus b\H,$  i.e. $[a] = \infty, [b] = 0,$ and $[c] = 1.$  Then any other point $p \in \HP^1$ may be written in coordinates 
\begin{equation}
p=[a\lambda + b]
\end{equation}
 with $\lambda \in \H \cong [a\H + b] = \HP^1 \setminus \{\infty\}.$   Now, as $c$ was chosen freely,  consider $c$ scaled by $c \mapsto \tilde c = cq \in \gamma$ for $q \in \H,$ and define a new normalized basis $\{\tilde a, \tilde b\}$  for $\H^2.$  It is clear that $\{\infty, 0, 1\}$ are preserved by this scaling. However,  $p$ is now given 
\begin{equation}\label{coor}
p = [\tilde a q^{-1}\lambda q + \tilde b].
\end{equation} 
 so that 
\begin{equation}\label{nocr}
\lambda \mapsto q^{-1}\lambda q
\end{equation}
under the scaling action. Thus we see that unlike $\CP^1,$ as a consequence of the noncommutativity of the quaternions,  there is no unique affine chart defined by a choice of three points.  Instead there is a family of projective coordinate systems for $\HP^1$ that assign those three points  to $\{\infty, 1,0\}.$  

Let $\{p_1,p_2,p_3\}$ be points in $\HP^1$ and  $\lambda \in \H.$ Then, we write $\lambda = x + n y$ where $x,y \in \R \subset \H$ and $n$ is a unit imaginary quaternion.   Thus, from (\ref{coor}), a  fourth point has coordinate $\lambda = x +  q^{-1} n q y,$  where $q$ is the free scaling parameter attached to the choice of coordinate system.  One moves through the family of coordinate systems by allowing $q$ to range over the quaternions, acting on $n \in Im( \H)$ by the representation of $SU_2$ on $\Im \H \equiv \R^3.$  Therefore the family of coordinate systems is equivalent to the orbit of $n$ under this action: the unit two-sphere.  One might define $\lambda$ to be the quaternionic cross-ratio in analogy with $\CP^1.$  However, as seen in its construction, $\lambda$ is not invariant under M\"obius transformations of $\HP^1.$   It is the quantities $x$ and $y$ which are invariant. Thus, in \cite{bobpink} this pair of real numbers is defined to be the quaternionic cross-ratio of four points in $\HP^1.$ However, the M\"obius invariance of this cross-ratio is obtained at the cost of uniqueness: specifying a cross-ratio and three points does not define a fourth point uniquely.
          
Given the coordinates of four points with respect to a coordinate system $\{q_1, q_2, q_3, q_4 \}$ in $\H \subset \HP^1,$  it is possible to  compute their \emph{quaternionic cross-ratio} writing  $[q_1, q_2, q_3, q_4 ] = (q_1 -q_2)(q_2 -q_3)^{-1}(q_3 - q_4)(q_4 - q_1)^{-1}.$ However,  this quaternionic cross-ratio is not a M\"obius invariant of $S^4:$ 
\begin{Pro}
Let the coordinates of four points in $ \HP^1$ be given with respect to a coordinate system by $\{q_1, q_2, q_3, q_4 \}$ in $\H$ so that their cross-ratio is computed by
\begin{equation} \label{crformula}
 [q_1, q_2, q_3, q_4 ] = (q_1 -q_2)(q_2 -q_3)^{-1}(q_3 - q_4)(q_4 - q_1)^{-1} = \lambda \in \H.
\end{equation}  
Let $\Phi$ be a M\"obius transformation of $\HP^1,$ then 
\begin{equation} \label{crformula}
 [\Phi(q_1), \Phi(q_2), \Phi(q_3),\Phi(q_4) ] = q \lambda q^{-1}
\end{equation}  for some $q \in \H$
\end{Pro}
{\em Proof}: The action of the M\"obius transformation $\Phi$ and then computation of the cross ratio is equivalent to a change of coordinates on $\HP^1$ preserving $\{\infty, 1, 0\}$ so that \ref{nocr} gives the result.$\square$

\begin{Cor}
Suppose that the cross ratio of four points $[q_1, q_2, q_3, q_4 ]  = \lambda \in \R.$ Then, $\{Re(\lambda), |Im(\lambda)|\}$ is invariant under M\"obius transformations of $\HP^1.$
\end{Cor} Thus, the quaternionic cross-ratio of four points is only well-defined if for those four points it is a real number.

Now,  choose an affine chart so that four points have coordinates $\{q_1,q_2,q_3,q_4\}$ and assuming $q_4 = \lambda = x + \hat n y \notin \R,$ one obtains that $\{q_1,q_2,q_3,q_4\} \subset \Span\{1,n\} \subset \H$ where $x,y \in \R$ and $n$ is a unit imaginary quaternion.
Thus four points in $\HP^1$ determine uniquely the two-sphere $\Span\{1,n\}\cup \{\infty\} \subset \HP^1$ provided that $y \neq 0$ i.e. there is no circle containing all of them. We may then choose coordinates so that $\lambda = x + i y \in \C \subset \H $ as it is always possible to find a $q \in \H$ such that $q^{-1} \hat n q = i.$   We see then that given four points in $\C \subset \H,$  it is always possible to choose coordinates so that the quaternionic cross-ratio is a complex number.

Assume $\lambda = x + iy \notin \R$ and consider the point \begin{equation} \label{4coor} \lambda_q = q \lambda q^{-1} = x + q i q^{-1} y.\end{equation} The real span of $\{1,\lambda_q\}$ is a  two-dimensional plane in $\H \cong \R^4$ which contains  the  ``real line" as the $\R-$span of $1.$  However, scaling the projective coordinates for $\HP^1$ by $q \in \H$ preserves $\{\infty,1,0\}$  and takes $x + q i q^{-1} y \mapsto \lambda.$ Thus, we associate the cross-ratio $\lambda$ to the point $\lambda_q$ and to a single complex number $\lambda$ there is a real two-sphere  whose points  $\lambda_q$ in appropriate coordinates have  complex quaternionic cross-ratio $[\infty,1,0,\lambda_q] =\lambda$ in $\C = \Span_R\{1,\lambda\} \subset \H.$ Each point $\lambda_q$ defines a real two-plane containing the real line $\R \subset \C \subset \H$ and the two-sphere of points is just the orbit of $i \in \C \subset \HP^1$ under the action of $SU_2.$  

From a geometric perspective: three points in $\HP^1$ are contained in a unique circle. Identifying these points with $\{\infty, 0,1\},$ this circle corresponds, by M\"obius transformation, to the real line $\R \subset \H.$  A fourth point not in $\R$ now spans a real two-plane which may then be identified with $\C \subset \H$ by change of coordinates and thus has a complex cross-ratio.  Finally this two-plane may be recognized as a two-sphere containing $\infty$ in $\HP^1.$  Thus, the set of fourth points of complex quaternionic cross-ratio $\lambda$ with respect to a given three is identified with the set of two-spheres in $\HP^1$ each containing the circle through $\{\infty, 0,1\}.$  The set of two-spheres containing a given circle is a familiar construction from the previous section. In fact, given points $p,q,$ and $r$ in a two-sphere $S = [p\C + q] \cong \C,$ the cross-ratio $\lambda$ is the coordinate of a fourth point $[s]= [p \lambda + q].$  Thus, Corollary \ref{steiner} gives us: 
\begin{Pro}\label{cr}
 Let $\{[p\wedge pj],  [q\wedge qj],  [r\wedge  rj]\} \in S^4 \subset Q^4,$ and $\lambda \in \C.$ Suppose $S = [p \wedge q]$ corresponds to a two-sphere  which also contains the point $[r] \in \HP^1.$  Then   the Steiner cross-ratio of $\{ [p\wedge pj],  [q\wedge qj],  [r\wedge  rj], [(p\lambda + q) \wedge (pj \lambda + qj)]) \}$ is exactly $ \lambda,$  where $[(p\lambda + q) \wedge (pj \lambda + qj)]$ corresponds to a two-sphere half-touching $S$ at $\{[p\lambda + q], [p \bar \lambda + q]\}.$ 
\end{Pro}

\begin{figure}[h]
 \includegraphics{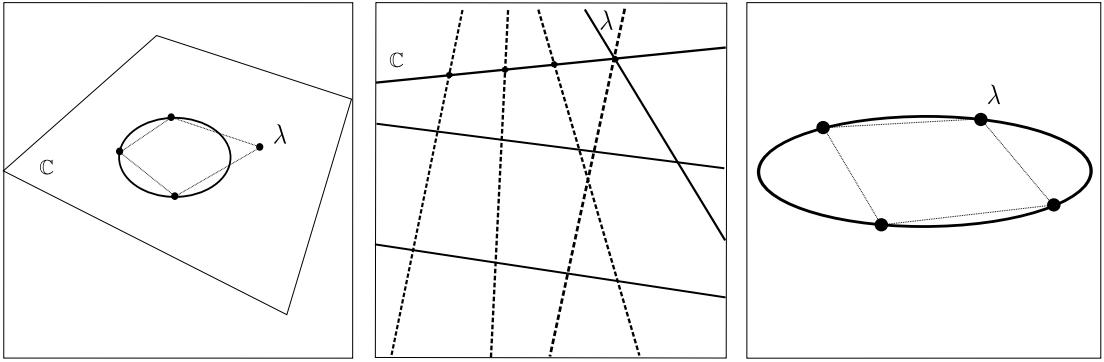} 

\caption{Four points with complex cross-ratio in $\HP^1,$ their twistor fibers in $\CP^3,$ and four associated points on a conic in $Q^4$ with the same cross-ratio.}
\end{figure}

We may now observe that, as a result of the noncommutativity of the quaternions,  quaternionic projective geometry  is fundamentally different from complex projective geometry.  We thus propose to ``complexify'' $\HP^1$ by considering $S^4$ as a real subset of $Q^4$ with respect to the real structure $\{j, Q^4\}.$ Circles in $S^4$ become nondegenerate planar conics in $Q^4.$ Then, given three points in $S^4 \subset Q^4,$ a choice of $\lambda \in \C$  uniquely specifies a fourth point in $Q^4$ whose Steiner cross-ratio, with respect to  the regulus in $Q^4,$ is $\lambda.$  If $\lambda \notin \R$ then, in M\"obius geometric terms, this  ``fourth point'' in $Q^4$ corresponds to a two-sphere  $S_\lambda \subset S^4.$  Choosing an oriented two-sphere $S$ containing the circle in $S^4$ defined by $\{p,q,r\}$ (incident to the regulus of $\{p,q,r\}$) and $\lambda \in \C$ defines a unique point in $S^4 \subset Q^4$ as the right-touching point in $S \cap S_\lambda.$ Thus the Steiner cross-ratio has the uniqueness property for circular points in $S^4 \subset Q^4.$  However, by ``complexifying'' $\HP^1$ we extend M\"obius geometry to twistor geometry and lose the second property of the complex cross-ratio: three points in $S^4 \subset Q^4$ do not define a point (M\"obius) transformation in $S^4$ but rather a general transformation of $Q^4.$

\subsection{Twistor geometry is the complexification of Lie geometry.}\label{3dim}

Before proceeding to develop discrete integrable systems in twistor geometry it will be useful to show how twistor geometry generalizes Lie geometry in $S^3.$   It is clear that if we restrict ourselves to a fixed $S^3 \subset S^4,$ then we may explore three-dimensional M\"obius geometry in $S^3$ by those M\"obius transformations of $S^4$ which preserve the given $S^3.$  The set of two-spheres and points in $S^3,$ whcih will be represented by the Lie quadric,  is clearly a subset of the set of two-spheres (and points) in $S^4.$  We will define Lie quadric as the real set of a real structure on  $Q^4.$  Lie geometry is an extension  of three-dimensional M\"obius geometry to include transformations taking points in two-spheres. From the perspective of Lie geometry, M\"obius geometry is given by ``point'' transformations, i.e. symmetries of the Lie quadric which take points in $S^3$ to points.  The space of null-lines in the Lie quadric is an example of a contact geometry. These contact lines will be considered as the real part of complex null lines in $Q^4.$  The geometry of $Q^4,$ that is of two-spheres in $S^4,$ is thus obtained by complexifying Lie geometry. However, as we have seen, the contact geometry of two-spheres in $S^4$ is not simply the complexification of the contact lines in the Lie quadric. Finally, we show that the three-dimensional complex quadric $Q^3,$ employed by Bryant in \cite{bry} and interpreted as the space of circles in $S^3,$ is understood naturally in terms of the full twistor geometry of $Q^4.$

 We will now indicate how to include three-dimensional geometry in $\HP^1.$ Given  an affine chart $\H \subset \HP^1,$ it is natural to identify $S^3$ with $Im(\H) \cup \{\infty\},$ where $Im(\H) = \Span_\R(i,j,k)\cong \R^3.$  Then, in analogy to the theory of circles in $S^2 \cong \C \cup \{\infty\}$ \cite{cir}, we may consider the nondegenerate quaternionic hermitian form $\mathfrak{h}$ defined by
\begin{equation*}
\begin{pmatrix}
0 && 1 \\ 1 && 0
\end{pmatrix}
\end{equation*}
with the properties that $\mathfrak{h}(v,w) = \overline{\mathfrak{h}(w,v)} \in \H$ and $\mathfrak{h}(v \lambda,w\mu) = \bar{\lambda} \mathfrak{h}(v, w) \mu$ for $v,w \in \H^2$ and $\lambda,\mu \in \H.$ Then, we compute the null-points of $\mathfrak{h}$ in $\H$ by 
 \begin{equation}
\begin{pmatrix}
\bar{z} && 1
\end{pmatrix}
\begin{pmatrix}
0 && 1 \\ 1 && 0
\end{pmatrix}
\begin{pmatrix}
z \\ 1
\end{pmatrix}
 = \bar{z} + z = 0.
\end{equation}
Thus, with respect to these coordinates up to a choice of real scale, we identify $S^3 = Im(\H) \cup \{\infty\}$ with $\mathfrak{h}.$

Now, the quaternionic hermitian form $\mathfrak{h}$ takes values in $\H = \C \oplus j\C.$ Thus, we write 
\begin{equation} 
\mathfrak{h}(v,w) = h(v,w)  + j \omega(v,w).
\end{equation}  Then, using $\mathfrak{h}(v,w) = \overline{\mathfrak{h}(w,v)},$ a short calculation shows that

\begin{equation}
\begin{split}
h(v,w) & =  \overline{h(w,v)} \\
\omega(w,v) & = - \omega(v,w).
\end{split}
\end{equation} 

The following may then be obtained by elementary computation \cite{mob}:

\begin{Lem} \label{frakh}
 Let $\mathfrak{h} = h + j \omega$ be a nondegenerate quaternionic hermitian form on $\H^2 \cong \C^4,$ then
\begin{enumerate}
 \item $h$ is a signature $(2,2)$ hermitian form on $\C^4.$
\item $\omega \in \bigwedge^2 {\C^4}^* $
\item $\omega(x,y) = -\overline{h(x,yj)}, $ and $ \omega(x, yj) = \overline{h(x,y)}. $
\item $h(x,xj) = 0,$ and $h(x,x) = \overline{\omega(x,xj)}$

\end{enumerate}
 
\end{Lem}

It is clear that a projective line $l \subset \CP^3$ corresponds either to a point in $S^3$ or to a two-sphere contained in the given three-sphere if and only if  for each $[x] \in l,$ $\mathfrak{h}(xq) = 0$ for $q \in \H.$  Thus, for each $[xq] \in [x \wedge xj] \subset \CP^3,$ $h(xq,xq) =0.$  Then, using Lemma \ref{frakh}, $l$ corresponds to a two-sphere contained in the three-sphere associated to a quaternionic hermitian form $\mathfrak{h} = h + j \omega$ if and only if  $h(x,x)=0$ for each $[x] \in l.$  Now, the hermitian perpedicular space $l^\perp = \P (\{x \in \C^4$ such that $h(x,l) = 0\})$ is a projective line.   Thus it is clear that $l = l^\perp$ if and only if $h(x,x) = 0$ for all $[x] \in l.$    Hence, we obtain
\begin{Pro}
Let $l$ be a projective line in $\CP^3.$ Then $l$ corresponds to a two-sphere (or point) contained in the three-sphere if and only if $l = l^\perp$ with respect to the hermitian form $h$ on $\C^4.$
\end{Pro}

Thus if we consider $ l \mapsto l^\perp$ as a map from $Q^4 \to Q^4,$ then the fixed points of that map correspond to two-spheres in the three sphere, including points of $S^3$ as two-spheres with zero radius.  We claim that this map defines a real structure on $Q^4$ so that the set of spheres in $S^3$ is given as the real set corresponding to that real structure. 

\begin{The}
Let $l$ be a projective line in $\CP^3,$ then define $\rho: Q^4 \to Q^4$ by $\rho(l) = l^\perp.$  Then, $\{\rho, Q^4\}$ defines a real structure. The real set is a real quadric of signature $(2,4)$ and corresponds to the set of two-spheres (and points) in $S^3.$  
\end{The}

{\em Proof}:   The hermitian form $h$ on $\C^4$ defines an anti-linear map by $x \mapsto  h(x,-) \in {\C^4}^*$ which induces an anti-holomorphic map $\C P^3 \to {\C \P^3}^*.$   This induces a natural map from $Q^4$ to ${Q^4}^*,$ the set of projective lines in ${\C \P^3}^*,$ which we will denote by $l \mapsto l^\flat.$ Now, projective duality interchanges lines in $\CP^3$ and lines in ${\C \P^3}^*.$ Let $L \subset {\C \P^3}^*$ be a projective line, then the projective dual $L^*= \{[x] \in \CP^3$ such that $\lambda(x) = 0$ for all $\lambda \in L.$ Thus, if $l \subset \CP^3$ is a projective line, then ${l^\flat}^* = l^\perp = \rho(l).$  Now, projective duality is a holomorphic map  from $Q^4$ to it's dual in ${\wedge^2 \C^4}^*.$  It is naturally induced by the pairing between  ${\wedge^2 \C^4}^*$ and ${\wedge^2 \C^4}$ defined for $\alpha,\beta \in \wedge^2 \C^4$ by $(\alpha, \beta) \mapsto \alpha \wedge \beta \in \wedge^4 \C^4 \cong \C.$    Thus, $\rho$ is given as the composition of a holomorpic and anti-holomorphic map and ${l^\perp}^\perp = l.$   Hence, $\{\rho, Q^4\}$ defines a real structure whose real set consists of projective lines $l$ such that $l = l^\perp.$ 

Now consider $\{v,vj,w,wj\}$ as a basis for $\C^4$ such that $[v],[w] \in S^3 \subset \HP^1.$  With respect to this basis we consider $[v] = \infty \in S^3$ and $[w] = 0 \in \R^3 \subset S^3.$  Then, 
\begin{equation}\label{basis}
B= \{v\wedge vj, w\wedge wj, v \wedge wj, vj \wedge w, \frac{1}{2}(v \wedge w - vj \wedge wj ), \frac{1}{2}(v\wedge w + vj \wedge wj) \}
\end{equation}   
is a basis for $\bigwedge^2\C^4$ with the property that for $\alpha \in B,$ $\rho([\alpha])= [\alpha].$  Thus, we may lift $\rho$ to an anti-linear map $\tilde \rho : \wedge^2\C^4 \to \wedge^2\C^4$ such that $\tilde \rho(\alpha) =  \alpha.$ Then, with respect to this basis $\tilde \rho$ acts as complex conjugation on $\wedge^2 \C^4$ and the fixed points of $\rho$ are then given by the projective image of the real span of $B$ in $Q^4 \subset \P(\wedge^2\C^4).$  Thus, the real set has the induced quadric structure from $Q^4$ on $\wedge^2\R^4$ given as the real span of $B.$  An elementary computation obtains that the signature of this quadric is $(2,4). \square$

Let $x\H \in S^3 \subset \HP^1,$  then the twistor fiber of $x\H$ is the projective line $[x \wedge xj]$ and we know that $h(x,x) = 0.$ Now, recall that the choice of an $S^3 \subset S^4$ is equivalent to a quaternionic hermitian form $\mathfrak{h} = h + j\omega,$  where $\omega \in \wedge^2{\C^4}^*$ is an alternating 2-form.  Then, by Lemma \ref{frakh}, 
\begin{equation}
\omega(x\wedge xj) = \omega(x,xj) = \overline{h(x,x)} = 0.
\end{equation}     
Hence, the points of $S^3$ are constrained by the linear condition
$
\omega^o = \{ \alpha \in \wedge^2 \C^4$ such that $ \omega(\alpha) = 0\}
$
and, with respect to the previous coordinates,
\begin{equation}
\begin{split}
\omega^o = & \Span \{v\wedge vj +   w \wedge wj, v \wedge wj - w \wedge wj, v \wedge wj + vj \wedge w, v \wedge w, vj \wedge wj \}\\ = & (v \wedge wj - vj \wedge w)^\perp .
\end{split}
\end{equation}   
A brief computation now shows that the signature of the $Q^4$ quadric restricted to the real part of this linear subset is $(1,4).$ Thus, we have obtained
\begin{Pro}
Let $\mathfrak{h} = h + j\omega$ be a nondegenerate quaternionic hermitian form. Let $\omega^o = \{ \alpha \in \wedge^2 \C^4$ such that $ \omega(\alpha) = 0\},$  then with respect to the real structure on $Q^4$ induced by $h,$ the real part of $\P(\omega^o) \cap Q^4 $ determines a light-cone model of $S^3 \subset S^4 \subset Q^4.$
\end{Pro}

Now, let us examine the complex quadric $Q^3 = \P(\omega^o) \cap Q^4.$  As discussed, $Q^3$ contains $S^3$ as a real subset, thus we now consider points in $Q^3 \setminus S^3$.  We note first  that the right action by the quaternionic $j$ preserves   $\P(\omega^o) .$  Let $[\lambda] \in \P(\omega^o) ,$  then  $\omega(\lambda j) = \overline \omega( \lambda) = 0.$  If $[\lambda] \in  Q^3 \setminus S^3$ then we recall that $[\lambda j]$ corresponds to the same two-sphere in $S^4$ with opposite conformal structure.  Further,   $[\lambda] \in  Q^3 \setminus S^3$ implies that the projective lines corresponding to $[\lambda] $ and $[\lambda j]$ in $\CP^3$ must be disjoint.  Thus the projective line $\P(W)$ spanned by $W = \{\lambda, \lambda j\}$ in $\wedge^2 \C^4$ must intersect $Q^4$ transversely.  Hence,  $\P(W) \cap Q^3 = \{[\lambda], [\lambda j]\}.$ Now, given $\omega^o$ as the complex span of $V,$ as indicated previously, we may write $\lambda = X- iY$ where $X$ and $Y$ are linearly independent elements of the real span $\Span_\R(V).$ Then an elementary analysis of the light-cone metric induced on   $\Span_\R(V)$ shows that $\P(\Span_{\R}W^\perp  ) \cap S^3$ determines a circle $L$ contained in $S^3.$ Thus, we see that to each pair of points $\{[\lambda], [\lambda j]\} \subset Q^3 \setminus S^3$ there is associated a circle contained in $S^3.$ We may think of this pair as the two possible orientations of a circle in space.  Conversly, given a circle $L \subset S^3.$ Choose three points on $\{[\alpha], [\beta], [\gamma]\} \subset L,$ so that $W= \Span\{\alpha, \beta, \gamma\} \subset \omega^o.$  Then $\Span_\C(W^\perp) \cap Q^3 \setminus S^3 = \{[\lambda], [\lambda j]\}.$  Thus, we obtain
\begin{Pro}
 Let $L$ be a circle in $ S^3 \subset Q^3 = \P(\omega^o) \cap Q^4,$ then there is a unique pair $\{[\lambda],[\lambda j] \in Q^3 \setminus S^3\}$ with the property that the two-sphere in $S^4$ associated to $\{[\lambda],[\lambda j] \}$ contains $L.$ The set of oriented circles (and points) in $S^3$ is thus parametrized by the three-dimensional complex quadric $Q^3.$  
\end{Pro}

Thus we have the novel interpretation that the space of oriented circles in $S^3$ is parameterized by a $3$-parameter complex family of two-spheres in $S^4$ each intersecting $S^3 \subset S^4$ transversely at their corresponding circle.  The interested reader may compare the discussion of this material in \cite{fer}, \cite{bry}, and \cite{bla}.

Now, consider the following classical result:

\begin{Lem}[Touching Coins Theorem]
Whenever four circles in three-space touch cyclically but do not lie in a common sphere, they intersect the sphere which passes through the four points of contact orthogonally. 
\end{Lem}

As the four circles do not lie on a common sphere, their representatives in $Q^3$ span a three-dimensional linear projective space.  Hence there is a unique two-sphere in $S^4$ in contact with all four representatives. With the generic assumption that the four points of contact do not lie on a circle, one obtains that they span the same three-dimensional space. Hence the unique two-sphere defined by those four points is also the unique two sphere in contact with all four representatives.  The two-sphere defined by the four contact points must  then half-touch each representative as the intersection with each contains exactly two points.  Thus, we interpret contact of circles in $S^3$ in terms of generalized contact of two-spheres in $S^4.$

\section{Discrete Integrable Systems in Twistor Space}

The theory of circular nets as a discrete integrable system in $\HP^1$ was developed in a paper by Bobenko and Pinkall \cite{bobpink} introducing the quaternionic cross-ratio and using a discrete Lax-pair formulation. Later, Hertrich-Jeromin, Hoffmann, and Pinkall \cite{pink} studied circular nets by the discrete evolution of a discrete curve using the quaternionic cross-ratio. This cross-ratio was assumed to be real and a further assumption determines the flow of the curve as  the discretization of an isothermic surfaces in $S^4.$ The idea is straightforward to relate.  

Starting with a discrete curve in $S^4,$ defined as an ordered chain of points $\{p_0,...,p_n\} \subset \HP^1,$ let $p_0^+ \in \HP^1$ and $\lambda_k \in \R$ for $ k\in \{0,...,n-1\}.$  Then, compute one step, $\{p_0^+,...,p_n^+\},$ of the discrete evolution of the initial curve by assuming that the quaternionic cross-ratio 
\begin{equation}\label{crsystem}
 [p_{k+1}, p_k, p_k^+, p_{k+1}^+] = \lambda_k
\end{equation}
is obtained.  We note that since $\lambda_k \in \R,$ the quaternionic cross-ratio is M\"obius invariant, uniquely determines  $p_{k+1}^+ \in \HP^1$ and the points $\{p_{k+1}, p_k, p_k^+, p_{k+1}^+\}$ are constrained to lie on circles. Specifying points $\{p_0^+, p_0^{++}...\}$ along a discrete curve transverse to the initial curve through $p_0$ determines a map from a subset of the square grid $\Z^2 \to \HP^1.$ With the assumption $\lambda_k = -1$  the image of this map is said to be ``discrete isothermic.''  The coordinate directions of $\Z^2$ define curvature lines on the discrete surface or surface patch.  For the sake of simplicity we will refer to the system defined by (\ref{crsystem}) where $\lambda_k$ has a constant value as the ``cross-ratio'' system. 

Hertrich-Jeromin et al., introduce the concept of Darboux transformation and Bianchi permutability for discrete isothermic surfaces \cite{pink}. Having obtained a discrete surface as the flow of a discrete curve, they obtain a Darboux transformation of the surface in a transverse direction to the surface by extending the cross-ratio condition in the direction of the transform.  The cross-ratio condition is seen as a discrete Riccati type equation. The permutability of Darboux transforms of discrete surfaces is seen as a symptom of ``integrability.'' 

We will now develop the theory of conjugate or ``quadrilateral'' nets in $Q^4,$ imposing a reality condition to obtain nets in $S^4.$  In particular, we will show that circular nets in $S^4$ are a ``real reduction'' of a generalized cross-ratio system which is given as a conjugate net in $Q^4.$ We will obtain as a consequence of this that the {\em complex} cross-ratio system in $S^2 \cong \CP^1$ is given by a conjugate net.  In addition, we will generalize the theory of principal contact element nets in the Lie quadric, developed in \cite{bob2}, to half-contact in $S^4,$ that is, generalized complex contact elements in $Q^4.$  

We first summarize the elements of the theory of conjugate nets as developed by A. Doliwa and A. Bobenko and show that discrete integrability is compatible with the addition of a real structure as a constraint.   

\subsection{Conjugate Nets}
Discrete conjugate nets are a basic element in recent theories of discrete differential geometry \cite{bob1}. Their introduction goes back to Sauer and they were developed extensively by Doliwa. They may be described as solutions of a set of partial {\em difference} equations. However, we prefer a purely geometric definition:
\begin{Def}
 A discrete conjugate net in $\CP^n$ is a  map $\phi : \Z^k \to \C^{n+1}$ with the property that the span of the image of the vertices on each 2-dimensional face in $\Z^k$ is three-dimensional. Equivalently, the image under $\phi$ of each face in $\CP^n$ is contained in a projective plane.    
\end{Def}

From the perspective of differential geometry,  $\Z^2$-nets parameterize discrete surfaces, with the directions of the square lattice coorresponding to the coordinate curves of the parameterization. Thus, fixing a direction and level surface, a $Z^3$-net consists of a discrete one-parameter family of surfaces. If we think of the level 1 surface as a transformation of the level 0 surface, one obtains a characterization of the fundamental transformation theory of discrete surfaces: two surfaces are fundamental transformations of one another if the resulting $\Z^2\times\{0,1\}$-net maintains the underlying properties of the two  \cite{dol2}.  Thus, two discrete conjugate nets are fundamental transforms of one another if the resulting net is a discrete conjugate net.   

Fundamental transformations of discrete conjugate nets exhibit Bianchi-like permutability, a fundamental property of  many integrable Hamiltonian PDEs related to classical surface theory \cite{bob1} e.g.  Darboux transforms of isothermic surfaces, B\"acklund transforms of $K$-surfaces.  Given a discrete conjugate net and two fundamental transforms of it, there exists a fourth conjugate net which is the fundamental transform of the two previous. The combinatorics of the permutability of solutions to the discrete conjugate net system is built into the elementary combinatorics of the $\Z^2$ and then $Z^n$ lattice. 

The basic lemma for Bianchi permutability in the discrete case is the following \cite{dol4}:
\begin{Lem}[Hexahedron Lemma]\label{hexa}
Given seven points $[\phi], [\phi_i], [\phi_{ij}]$ ($1 \leq i < j \leq 3$) in $\CP^n,$ such that each of the three quadrilaterals $\{ [\phi], [\phi_i], [\phi_{ij}],[\phi_j]\}$ is planar, define three projective planes $\tau_i\Pi_{jk}$ as those passing through the point triples $\{[\phi_i], [\phi_{ij}],[\phi_{ik}]\}$, respectively. Then these three planes intersect at one point:
\begin{equation}
 [\phi_{123}] = \tau_1\Pi_{23}\cap \tau_2\Pi_{13} \cap \tau_3\Pi_{12}
\end{equation}
\end{Lem}

This result is readily obtained for conjugate nets \cite{dol4} as the intersection of three planes in space is generically a point. However for discrete circular nets the equivalent Lemma for seven circles in classical geometry is refered to as ``Miquel's Lemma'' \cite{dol5}, \cite{ped}. The crucial step in obtaining Bianchi permutability is that the Hexahedron Lemma
may be extended from the $\Z^3$ lattice to a $\Z^4$ lattice.

\begin{figure}[h]
 \includegraphics{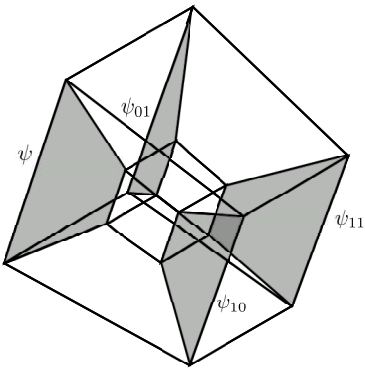}
 \caption{Bianchi permutability of elementary faces of four discrete conjugate nets in the elementary hypercube of  the $Z^4$-lattice. }
\end{figure}

\begin{Pro}[Bianchi Permutability]
Let $\psi$ be a discrete conjugate net. Suppose that there exist discrete conjugate nets $\psi_{01}$ and $\psi_{10}$ which are each fundamental transformations of $\psi.$ 
Then, there is a two-parameter family of planar discrete conjugate nets $\psi_{11}$ with the property that $\psi_{11}$ is a fundamental transformation of $\psi_{01}$ and $\psi_{10}.$  
\end{Pro}

{\em Proof}: One may arrange the nets $\phi, \xi, \eta$ so that an elementary face of $\phi$ is related to an elementary face of $\xi$ and $\eta$ along independent directions in the four-dimensional lattice with $\tau_1 \phi = \eta$ and $\tau_2 \phi = \xi$. Now, consider the plane spanned by the points $\phi, \xi, \eta.$ There is a two-parameter choice of points in this plane which will be labeled $\psi.$ By the hexahedron lemma there is now a unique point $\psi_{123}$ so that $\{\phi_3, \xi_3, \eta_3,\psi_{123}\}$ is planar. Now, by the hexahedron lemma again there is a unique point $\psi_{124}$ such that $\{\phi_4, \xi_4, \eta_4,\psi_{124}\}$ is planar. The proof then follows that of the four-dimensional hexahedron lemma defining a final point $\psi_{1234}$ and determining an elementary face of the planar net $\psi. \square$

Discrete circular nets in the light-cone model of $S^4$ are quadratic reductions of conjugate nets \cite{dol3}.   Assume that the vertices of a conjuate net are mapped into the light-cone in $\R^6$, the zero set of the Lorentz quadratic form. Then, requiring that the quadratic form restricted to each face of the net be irreducible, the intersection of each face-plane of the conjugate net with $S^4 \subset \RP^5$ is a circle.
What is then required is the compatibility of this quadratic constraint with the Hexahedron Lemma: the eighth point must lie on the quadric.  This compatibility, for quadrics in general,  comes as a result of the classical ``eight points lemma'' \cite{ped}, where opposing ``parallel'' sides of the hexahedron are considered as degenerate quadrics. Thus, discrete integrability for circular nets is obtained as a consequence of applying a quadratic constraint to an underlying conjugate net. 

We now consider $S^4$ as a real quadric variety contained in the complex quadric $Q^4 \subset \CP^5.$ Additionally, $S^4$ is considered as the real set corresponding to the real structure induced on $\bigwedge^2\C^4$ by the right action of $j\in \H$ on $\C^4 \cong \H^2.$  Thus, we require that the real structure be compatible with the Hexahedron lemma on the underlying complex quadric.  


\begin{Lem}
Consider seven points $[\phi], [\phi_i], [\phi_{ij}]$ ($1 \leq i < j \leq 3$) in a complex quadric hypersurface $Q^{n-1} \subset \CP^n$ with a given real structure and non-empty real set.  Suppose that these points are contained in the real set. Then, under the conditions of Lemma \ref{hexa}, there is a unique eigth point contained in the real set. 
 \end{Lem}
{\em Proof}:
Under the conditions of Lemma \ref{hexa} each of the seven points are coplanar with their lattice-neighbors. Hence  they span a three-dimensional projective space and it suffices to consider $n=3.$ The result may then be obtained by direct computation using an explicit formula in coordinates given by a basis of representatives of the intial four points of the  lattice hexahedron: $\{\phi_1, \phi_2, \phi_3,\phi_3\}.$  These points are invariant with respect to the real structure so that  the initial seven points have real coordinates. The eighth point is then given by $\phi_{123} = y_0  \phi+ y_1 \phi_1 + y_2 \phi_2+ y_3 \phi_3$ for 
\begin{equation}\label{doliwa}
 \begin{split}
y_0 & = a_0 b_0 c_0 ( \frac{1}{a_2 b_1 c_3} + \frac{1}{a_1 b_3 c_2}), \\ 
y_1 & = \frac{b_0 c_0}{b_3 c_2} +\frac{a_0 c_0}{a_2 c_3} +\frac{c_0^2}{c_2 c_3}, \\
y_2 & = \frac{a_0 b_0}{a_1 b_3} +\frac{b_0 c_0}{b_1 c_3} +\frac{b_0^2}{b_1 b_3}, \\
y_3 & = \frac{a_0 c_0}{a_1 c_2} +\frac{a_0 b_0}{a_2 b_1} +\frac{a_0^2}{a_1 a_2}.
 \end{split}
\end{equation}
where
\begin{equation}
 \phi_{12} = \begin{bmatrix}
         a_0 \\ a_1 \\ a_2 \\ 0
        \end{bmatrix},
\phi_{13} = \begin{bmatrix}
         b_0 \\ b_1 \\ 0 \\ b_3
        \end{bmatrix},
\phi_{23} = \begin{bmatrix}
         c_0 \\ 0 \\ c_2 \\ c_3
        \end{bmatrix}.
\end{equation}
from \cite{dol1}. Conceptually, the projective plane corresponding to each face of the conjugate net contains a real projective plane given by the real structure so that the proof is obtained by satisfying the properties of the ``eight-point" theorem for real quadrics .  $\square$

Thus, we may construct discrete conjugate nets in $Q^4$ and impose a real structure as a constraint on the vertices which preserves the discrete integrability of the underlying conjugate net.  Doing this for the structure imposed by the quaternionic $j$ we obtain conjugate nets in $\RP^5 \subset \CP^5$ with vertices in the light-cone model for $S^4:$ discrete circular nets.  

Now, it will be useful to consider the further reduction of the discrete conjugate net system in $Q^4$ by restricting to a fixed two-sphere $S$ contained in $S^4 \subset Q^4.$ As $\HP^1 \cong S^4,$ a given two-sphere contained in $\HP^1$ has two twistor lifts, $S$ and $Sj \subset \CP^3$ and thus two representative points $[S],[Sj] \in Q^4.$ We identify the two-sphere $S$ as a subset of $S^4 \subset Q^4$ by determining all of the twistor fibers incident to the twistor lifts $\{S,Sj\}$ in $\CP^3.$ However, we first examine the set of all projective lines incident to $\{S,Sj\}.$ 
\begin{Lem}
The family of projective lines incident to the twistor lifts $S$ and $Sj \subset \CP^3$ is given by $Q^2_S = \P(\{S,Sj\}^{\perp})\cap Q^4,$ a $2$-dimensional complex nondegenerate subquadric of $Q^4 \subset \CP^5.$ 
\end{Lem}

The twistor fibers will now be determined by restricting the real structure $j$ from $Q^4$ to $Q^2_S.$  As the action on $Q^4$ by $j$ preserves $\{[S],[Sj]\},$ this restriction determines a well-defined real structure $\{j, Q^2_S\}$ with real points exactly those elements of the form $[a\wedge aj]$ with $[a] \in S \subset \CP^3.$  Thus, we identify the two-sphere $S$ as a subset of $S^4$ by the real set of this induced real structure.  

\begin{Lem}
 A fixed two-sphere $S \subset \HP^1$ is determined as a subset of $S^4 \subset Q^4$ by the real set of the real structure $\{j, Q^2_S\}$ where $j$ is the restriction of the right action of $j \in \H$ on $Q^4$ and $Q^2_S = \P(\{S,Sj\}^{\perp})\cap Q^4$ for some $S$ and $Sj$ homogeneous representatives of the points in $Q^4$ corresponding to the two conformal structures on $S.$  
\end{Lem}

The reduction of the discrete conjugate net system to $S^2 \subset S^4$ is obtained first by the linear constraint $\{S,Sj\}^{\perp}$ to the vertices of the net. This determines a discrete net in the complex quadric $Q^2_S$ where each face is contained in a conic.  If these conics are irreducible, the addition of the reality condition obtains a circle in $S^2 \subset Q^2_S$ in each face as the real set of the underlying conic, thus determining a discrete circular net in $S^2.$  It remains to interpret the general ``discrete conic net'' in $Q^2_S$ in terms of $S^2.$ 

\subsection{The Complex Cross-Ratio System in $S^2$ is given by a Discrete Conjugate Net}

The structure of $Q^2_S$ can be explained in terms of projective lines in $\CP^3.$   The quadric $Q^2_S$ is considered as the set of projective lines incident to the skew lines $\{S,Sj\},$ the two twistor lifts of the two-sphere $S \subset \HP^1.$ Thus, each point in $Q^2_S$ can then be written as $[a_1 \wedge a_2 j ],$ where $[a_1], [a_2] \in S \subset \CP^3$ and $[a_2 j] \in Sj.$ Let $[a] \in S$ then $[aj] \in Sj$ so that the twistor fiber $[a\wedge aj] \in Q^2_S.$ 
 
 Given the discrete curve $\{p_0,...,p_n\}$ with $p_k \in S \subset \HP^1,$ (\ref{crsystem}) defines the associated cross-ratio system for a given $\lambda \in \R.$ Now, let $p_0^+ \in S.$ We compute the unique point $p_1^+$ so that the cross-ratio condition is satisfied. As $\lambda \in \R$ and $\{p_1, p_0, p_0^+\}$ are twistor fibers, $p_1^+$ is a twistor fiber over a point in $S$  and the points $\{p_1, p_0, p_0^+, p_1^+\}$ lie on a circle.    

In $\CP^3,$ the three skew lines $\{[p_0^+ \wedge p_0^+j],[p_0 \wedge p_0j],[p_1 \wedge p_1j]\}$ are incident to $S$ and $Sj.$   They define a regulus containing the points $S$ and $Sj$ which sweeps out a complex quadric surface in $\CP^3.$ The regulus corresponds to the complex conic $C= \P(\Span\{[p_0^+ \wedge p_0^+j],[p_0 \wedge p_0j],[p_1 \wedge p_1j]\}^\perp) \cap Q^4 \subset Q^2_S.$ The line $S =[p_0 \wedge p_1] \subset \CP^3$ can be identified with $\C \subset \H \subset \HP^1$ so that the cross-ratio of points on the projective line in $\CP^3$ is the same as the the quaternionic cross-ratio restricted to $\C \subset \H.$ Now, Proposition \ref{cr} says that the Steiner cross-ratio with respect to $\{[p_0^+ \wedge p_0^+j],[p_0 \wedge p_0j],[p_1 \wedge p_1j]\}$ is exactly the cross-ratio with respect to the points $[p_0^+], [p_0], [p_1]$  in $[p_0 \wedge p_1].$ Referring back to (\ref{doliwa}) we see that as $\{[p_0^+ \wedge p_0^+j],[p_0 \wedge p_0j],[p_1 \wedge p_1j]\}$ are all fixed points of $j$ i.e. real points, the fourth point which satisfies the cross-ratio condition must also be real.  Thus, using the Steiner cross-ratio we have recovered the real cross-ratio system in $S^2.$  

Now,  suppose we choose $\lambda \in \C,$ so that $\lambda \notin \R \subset \C.$ The immediate consequence of choosing the cross-ratio to be a complex number is that, if we attempt to repeat the above computation, the point $p_1^+$ is no longer the twistor fiber of a point in $S$ but, corresponds to a two-sphere half-touching $S$ and $Sj.$  The twistor lift of $S$ in $\CP^3$ is in the regulus generated by $\{[p_0^+ \wedge p_0^+j],[p_0 \wedge p_0j],[p_1 \wedge p_1j]\}$ and the line corresponding to $p_1^+$ is in the dual regulus containing the generators. 
Then, by Corollary \ref{steiner}, $p_1^+$ corresponds to the unique point in $S,$ identified with $\C,$ having cross-ratio $\lambda$ with respect $\{p_0^+, p_0, p_1\}.$ However, we note that $p_1^+ \in Q_S^2$ is not a point $S.$ 
Thus, we have obtained a correspondence between the evolution of a discrete curve in $Q_S^2$ by complex Steiner cross-ratio and the evolution of a discrete curve in $S \cong \CP^1$ by complex cross-ratio. 
 
\begin{The} \label{complexcr}
Let $S$ be a two-sphere in $S^4$ and $Q^2_S = \P(\{S,Sj\}^{\perp})\cap Q^4.$ Then, with respect to a choice of conformal structure on $S,$ each discrete conic net in $Q^2_S$ corresponds to a complex cross-ratio net in $S.$ Given a discrete curve $\{p\} \subset S \subset Q^2_S$ and $\lambda \in \C,$ then there is a unique discrete conjugate net in $Q^2_S$ extending $\{p\}$ corresponding to the complex cross-ratio evolution of $\{p\}$ by $\lambda$ in $S.$  
\end{The}

The evolution of closed discrete curves in $\CP^1$  by the complex cross-ratio system was extensively investigated in \cite{per}. Computing  the ``holonomy'' around the discrete curve allowed the association of the cross-ratio $\lambda$ with a discrete ``spectral'' parameter.  Integrability is then defined  to be the invariance of the associated ``spectral curve'' over the evolution of the discrete curve. If we identify $\CP^1$ with the light-cone model of $S^2,$ faces with complex (non-real) cross-ratio do not lie within a projective plane in $\RP^3.$ However, as we identify $S^2$ with the real set of $Q_S^2,$ we may consider $Q_S^2$ as the complexification of the light-cone model of $S^2.$  What we have shown is that the complex cross-ratio system is determined by a discrete conjugate net in the complexified light-cone model. Thus, we have also shown that it satisfies the integrability condition of Bobenko and Doliwa.

\subsection{Generalizing the Real Cross-Ratio System in $S^4$}

We have shown that we can generalize the real cross-ratio system in $S^2$ by complexifying the associated discrete conjugate net system in the light-cone model of $S^2.$ The resulting discrete conjugate nets in the complex two-dimensional quadric $Q_S^2$ are completely associated with solutions of the cross-ratio system in $S^2$ with {\em complex} cross-ratio.  Thus, it is natural to consider generalizing the real cross-ratio system in $S^4$ by considering discrete conjugate nets in the complexification of the light cone model of $S^4.$ We show that  the relationship between circular nets and nets of contact elements, established by Bobenko and Suris \cite{bob2}, is generalized for these \emph{discrete conic nets} in $Q^4.$    

\begin{Def}[Discrete Conic Net]
A discrete conjugate net $\phi :\Z^k \to {Q} \subset \CP^n$ in a
quadric $Q$ such that the intersection of the planes of
elementary faces of the lattice with the quadric
are irreducible quadratic curves is called a discrete conic net.
\end{Def}

The immediate interpretation of a discrete conic net is that it is a net of spheres in $S^4.$ However, we have seen that for nets in $Q_S^2$ discrete conic nets are exactly identified with complex cross-ratio nets by choosing either their point of left or right contact with a fixed two-sphere $S.$  In the work of Bobenko and Suris in discrete differential geometry, nets of spheres play an essential role \cite{bob2}.  In particular, ``S-isothermic'' nets provide the foundation for the theory of discrete minimal surfaces \cite{bob3}.  Following the outline of 19th century surface theory, there is an intimate connection between contact geometry of spheres and special surfaces in space \cite{bur}. 
\begin{Def}[Contact Element] A contact element is the one-parameter family of two-spheres touching at a point. 
\end{Def} We have seen in Corollary \ref{touchingnull} that contact elements for $S^4$ are null lines in $Q^4$ with the property that each null line contains a real point, the representative of a twistor fiber.   

Given a smooth surface in space, it is natural to consider the map to the space of contact elements associating to each point on the surface the set of spheres touching at that point.  The intuition for the definition of principal contact element nets goes back to Blaschke \cite{bla}.  Travelling \emph{infinitesimally} in the principal curvature directions along the surface,  neighboring contact elements will have a common sphere, the principal curvature sphere \cite{bob1}.  Thus, we can consider a discrete contact element net as a map $\Phi$ from $\Z^2$ to the space of null lines in $Q^4$ such that $\Phi(m,n)$ is a contact element.  Each contact element contains a point in space as the real point under the action of quaternionic $j$ on $Q^4.$ However, in the smooth theory, the set of points taken from each contact element only defines a smooth surface if the Legendre condition is satisfied \cite{bur}.  Thus, in the discrete theory we define a principal contact element net, generalizing the definition in \cite{bob1}.

\begin{Def}[Principal Contact Element Net]
 Let $\Phi$ be a map from $\Z^2$ to the space of null lines in $Q^4$ such that $\Phi(m,n)$ is a contact element for $(m,n) \in \Z^2$ and each contact element contains a distinct real point. We say $\Phi$ is a principal contact element net if  $\Phi(m,n) \cap \Phi(m+1,n) \neq \emptyset  $ and  $\Phi(m,n) \cap \Phi(m,n+1) \neq \emptyset.$
\end{Def}

\begin{figure}[h] \label{phi}
 \includegraphics{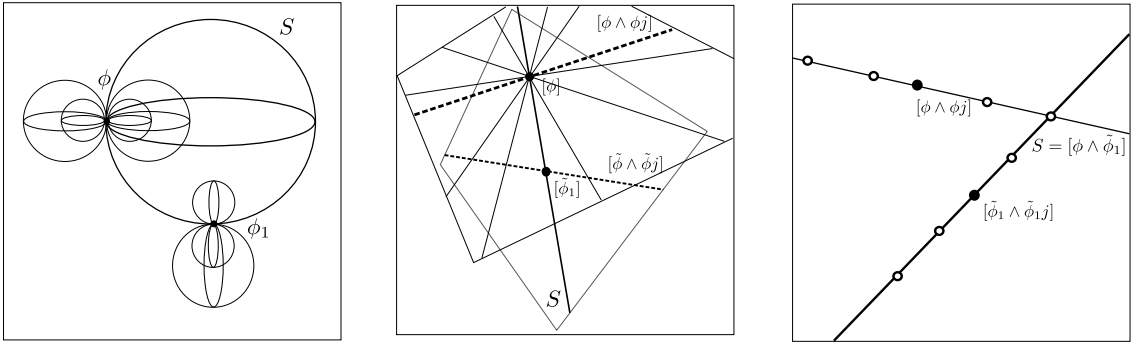}
 \caption{Choosing a contact element at $\phi$ induces a unique contact element at $\phi_1.$}
\end{figure}
Bobenko and Suris \cite{bob2} \cite{bob1} show that a principal contact element net in $S^3$ induces a circular net and the converse. We prove a generalization of this result in one direction, showing that circular nets induce principal contact element nets in $S^4.$  

\begin{The}\label{contact}
A principal contact element net is \emph{congruent} to a discrete net $\phi : \Z^2 \to \S^4$ if $\phi(k)$ is an element of the contact line at that lattice point.  Given a circular net $\phi: \Z^2 \to S^4,$ there exists a complex two-parameter family of principal contact element nets congruent to $\phi.$  
\end{The}
{\em Proof}: Let $\{\phi, \phi_1, \phi_2, \phi_{12}\}$ be an elementary face of a circular net in $S^4.$  From Theorem \ref{RNT}, a contact element at $\phi$ is defined by choosing a projective line $s$ incident to the twistor fiber of $[\phi \wedge \phi j] \subset  \CP^3.$  Assume that $s$ is chosen so that $s \cap [\phi_1 \wedge \phi_1 j] = \emptyset $ and  that $s \cap [\phi \wedge \phi j] = [\phi] \in \CP^3.$   This is equivalent to choosing a two-sphere in $s \subset S^4$ containing the point $\phi \in S^4$ and disjoint from the point $\phi_1.$ Then the contact element is the set of two-spheres touching $s$ at $\phi.$  This contact element is a null line $l \subset Q^4$ defined by the set of projective lines contained in the projective plane $\Pi$ in $\CP^3$ spanned by $\{[\phi \wedge \phi j], s\}$ and incident at $[\phi] \in \Pi.$  Following Corollary \ref{pointplane} we write $l = \{[\phi], \Pi\}.$ The twistor fiber of $\phi_1$ must intersect $\Pi$ at exact one point $[\tilde \phi_1] \in \CP^3.$ Then, the projective line $s_1$ spanned by $\{[\phi], [\tilde \phi_1]\}$ corresponds to a two-sphere in $S^4$ containing $\phi$ and $\phi_1$ as in Figure \ref{phi}.  Thus, there is a unique contact element $\{[\tilde \phi_1], \Pi_1\} = l_1 \subset Q^4$ at $[\phi_1]$ intersecting $l$ defined by the projective plane $\Pi_1$ spanned by $\{[\tilde \phi \wedge \tilde \phi j], S\}$ and the point $[\tilde \phi_1] \in \CP^3.$  We may similarly construct a contact element $l_2$ at $\phi_2.$  Note that the construction of contact elements above did not rely on the circularity of the face.

What remains to be shown is that this may be extended consistently around $\{\phi, \phi_1, \phi_2, \phi_{12}\}$  so that the contact lines induced at $\phi_{12}$ from $l_1$ and $l_2$ agree.  Now, the span of the lines $\{l,l_1,l_2\} \subset Q^4 \subset \CP^5$ is a three-dimensional projective space.   Then, by the circularity assumption, $\phi_{12}$ is contained in the projective plane spanned by $\{\phi, \phi_1, \phi_2\}$ in $Q^4$  which is contained in the span of  $\{l,l_1,l_2\} \cong \CP^3.$ Hence, the span of $\{\phi_{12}, l_1\} \subset  \CP^3$ must intersect the line $l_2 \subset \CP^3$ in exactly one point $s_{12} \in Q^4.$ But, the  span of $\{\phi_{12}, l_1\}$ is induced by the set of lines intersecting at $\phi_{12}$ and incident to $l_1,$ the set of contact lines at $\phi_{12}.$ Thus, there is a unique contact line $l_{12}$ spanned by $\{\phi_{12}, s_{12}\}$ intersecting $l_1$ and $l_2.$  

Finally, there is a two parameter of contact elements $l$ defined as $\{$point, plane$\}$ pairs  along the twistor fiber of $\phi$ and the result is obtained. $\square$   

Principal contact element net consist of lines in the space of two-spheres in $S^3:$  the Lie quadric.  We define the Lie quadric in Section \ref{3dim} as the real set of the real structure on $Q^4$ induced by a choice of $S^3 \subset S^4.$ Given a principal contact element net, each contact element is a real projective line parameterizing the set of two-spheres in $S^3$ tangent at a point. Thus, each contact element contains a point in $S^3$ and these points form a discrete net $\phi: \Z^2 \to S^3.$  Bobenko and Suris prove the converse of   Proposition \ref{contact}: given a principal contact element net in $S^3$, the associated discrete net is circular.  However, this is not true for principal contact element nets in $S^4.$  In fact this may be shown for a complex-cross ratio system in a fixed two-sphere $S \subset S^4.$

\begin{figure}[h]
\includegraphics{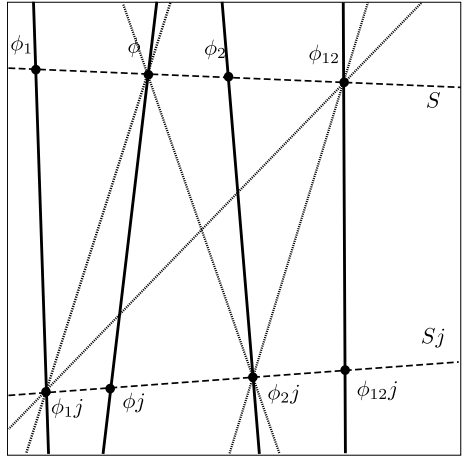}
\caption{Inducing a principal contact element net on a non-circular net.}
\end{figure}
\begin{Pro}
Given a two-sphere $ S \subset S^4,$ let $\phi: \Z^2 \to S$ be a discrete complex cross-ratio net. Then, there exists a principal contact element net in $S^4$ congruent to $\phi.$ 
\end{Pro}

{\em Proof}:  Let $\{\phi, \phi_1, \phi_2, \phi_{12}\}$ be an elementary face of a cross-ratio net in $S.$  Recall that a contact element in $S^4$ is defined by a $\{point, plane\}$ pair in $\CP^3.$ Define the initial contact element at $\phi$ by $l = \{ [\phi ], \Pi\},$ where $\Pi$ is the projective plane spanned by the lines $\{ [\phi \wedge \phi j], Sj\}.$  Following the construction of Prop. \ref{contact}, the twistor fiber of $\phi_1$ intersects $\Pi$ at $[\phi_1 j].$ Let  $\sigma_1$ be the two-sphere represented by the line in $\CP^3$ spanned by $\{[\phi], [\phi_1 j] \}.$  Thus, $\sigma_1$ is 
the unique two-sphere right-touching $S$ at $[\phi]$ and $Sj$ at $[\phi_1 j].$ Note that $\Pi$ is also the span of $\{ [\phi \wedge \phi j], \sigma_1\}.$ 
Then, $\sigma_1$ defines the contact element $\{[\phi j], \Pi_1\},$ where $\Pi_1$ is spanned by $\{ [\phi_1 \wedge \phi_1 j], \sigma_1\}.$  Thus, by construction $l \cap l_1 = \{\sigma_1\}.$ Notice that $\Pi_1$ is also the span of  $\{ [\phi_1 \wedge \phi_1 j], S\}.$   Now, repeat this construction for the second point so that $l_2 $ is defined at the point $ [\phi_2 j]$ by the plane $\Pi_2 = \{  [\phi_2 \wedge \phi_2 j], S\}$ and $\sigma_2 = [\phi \wedge \phi_2 j] .$  Finally, $[\phi_{12} \wedge \phi_{12}j] \cap \Pi_1 = [\phi_{12} \wedge \phi_{12}j] \cap \Pi_2 = \{[\phi_{12}\}$ so that the construction induces a unique contact element $l_{12} = \{ [\phi_{12}], Sj\}.$   We observe that $l_{12} \cap l_1 = {\sigma_1}^2  = [\phi_{12} \wedge \phi_1 j]$   and $l_{12} \cap l_2 = {\sigma_2}^1  = [\phi_{12} \wedge \phi_2 j].$  Now, observe given an initial choice of contact element $\{[\phi], \{[\phi]\wedge [\phi j], Sj\}$ or $\{[\phi j], \{[phi]\wedge [\phi j], S\},$  the principal contact element net is determined by the net $\phi.$ $\square$      
  

This is not a counterexample to the original result of Bobenko and Suris because we have been able to construct a principal contact element net where each contact element lies outside of any three-sphere containing $S.$ The contact elements in fact lie in the quadric $Q_S^2,$ so that each is defined by the sphere touching a two-sphere half-touching $S$ and $Sj.$  Thus, this construction depends upon the geometry of two-spheres in $S^4.$  As discussed previously, the cross-ratio net in $S$ is given by a conjugate net in $Q_S^2.$  Thus we have shown that it always possible to construct a principal contact element net in $S^4$ to a cross-ratio net in a fixed two-sphere.

\section{Conclusion}

The complex cross-ratio system in $\CP^1$ has been developed in \cite{per} as an algebraic integrable system. For solutions periodic in at least one direction in $\Z^2,$  there is a spectral curve assigned via the holonomy problem.  Described by (\ref{crsystem}), the cross-ratio $\lambda$ becomes the spectral parameter on the spectral curve.  Thus, formulating discrete curve evolution as an algebraic integrable system, the transition from real to complex cross-ratio is  equivalent to deforming the spectral parameter of the underlying system from real to complex values.  

The cross-ratio system in $\CP^1$ is defined for both real and complex values of the cross-ratio.  However,  identifying $ \CP^1$ with $S^2,$ for real values of the cross-ratio the cross-ratio system is given in purely geometric terms by circular nets, conjugate nets in the light-cone model of $S^2.$     Thus, it would be natural to attempt to encompass complex values of the cross-ratio by complexifying the light-cone model of $S^2$ and then considering complex conjugate nets in the complex light-cone as a two-dimensional complex quadric.  But,  these complex conjugate nets generally consist of points which are non-real and do not lie on $S^2.$  We have shown that the complexification of the light-cone model of $S^2$  is naturally defined in terms of twistor geometry as the quadric $Q_S^2.$   In Theorem \ref{complexcr} we define a natural transformation between conic nets, as solutions of the conjugate net equations in $\CP^4$ constrained to $Q_S^2$ and solutions of the cross-ratio system in $\CP^1$ for complex and real values of the cross-ratio.   This confirms the approach of Doliwa in \cite{dol1} and provides a concrete example of a discrete integrable net in a Grassmannian \cite{bob4}.  We claim that this conic net system in the complexified light-cone should be related in a natural way to the algebraic discrete integrable system in \cite{per}.  

Now, the conic net system is a reduction of a `master system' defined by conic nets in $Q^4,$ the space of two-spheres in $S^4.$  We see that ``point-lattices,'' lattices made up of points in $S^4$ (or $S^2 \subset S^4$,) are identified by the reality condition determined by the right-action of the quaternion $j$ on projective lines in $\CP^3.$  We have shown that this master system generalizes the ``complex'' quaternionic cross-ratio system in $S^4$ defined by Bobenko and Pinkall in \cite{bobpink} for {\em real} values of the cross-ratio.  Solutions for real values of the cross-ratio correspond to ``curvature'' parameterized discrete surfaces and ultimately discrete isothermic surfaces in $S^4.$ We determine that, in general, solutions for complex values of the cross-ratio correspond to discrete surfaces made up of spheres in $S^4.$  This suggests that deformations of the spectral parameter from real to complex, in the smooth case for conformal evolution of surfaces in $S^4$ \cite{sch}, determine surfaces made out of spheres.  

Bobenko, Hoffmann, and Springborn in their paper ``Minimal Surfaces from Circles,'' \cite{bob3}  define discrete minimal surfaces in $\R^3$ using a discrete net and  dual discrete net, where the image of the map at a vertex of each lattice is either a circle or sphere in $S^3.$  In developing three-dimensional conformal geometry as a subgeometry of four-dimensional conformal geometry, we determine that the space of circles in $S^3$ is naturally identified with a subset of $Q^4:$ two-spheres in $S^4.$  Thus we note that, from the perspective of twistor geometry,  the nets and dual nets defining ``$S$-isothermic nets'' each specify points into $Q^4$ either corresponding to a two-sphere contained in $S^3$ or a two-sphere in $S^4$ which interesects $S^3$ in a circle.  This suggest there may be a more general definition of $S$-isothermic nets in $S^3$ given in terms of four-dimensional geometry.  A natural direction to explore would be to attempt to deduce a geometric interpretation for the generalized isothermic nets of Doliwa \cite{dol1}  defined in $Q^4.$ These generalized isothermic nets would be a specialization of the conic net system we have described.            

Implicit in our discussion of three-dimensional geometry is a classical result known as Lie's ``line-sphere correspondence.''  This is an association between lines in $\RP^3$ and two-spheres in $S^3.$  The real structure and associated invariant basis (\ref{basis}) may be used as a basis for $\wedge^2\R^4 \subset \wedge^2 \C^4$ and to define a {\em real} Pl\"ucker quadric contained within the complex Pl\"ucker quadric.  The correspondence is defined explicitly by choosing affine charts to describe the intersection of the real Pl\"ucker quadric with the Lie quadric, given as real set. This provides a dictionary relating results about ``asymptotic'' surfaces in $\R^3 \subset \RP^3$ and surfaces in Lie sphere geometry.   Thus, the paper of Doliwa \cite{dol6} may be extended and related to work on S-isothermic surfaces in Lie geometry \cite{bob2}.   

Konopelchenko and Schief have published several papers on the discrete (quaternionic) KP equations \cite{ks1} \cite{ks2} \cite{ks3}.   The discrete KP hiearchy is alternately defined in terms of a quaternionic multi-ratio condition or the classical Clifford point-circle configuration \cite{ks3}. Twistor geometry, as we have developed it in this paper, should provide a natural framework for interpreting and extending these results.  

Finally, it remains to fully determine the relationship between the conic nets we have defined and the algebraic integrable system determining the cross-ratio system in $S^2.$  This should allow us to extend the analysis in \cite{per} to discrete surfaces in $S^4$ and provide a tool for investigating the integrable evolution of smooth surfaces in $S^4.$

\thanks{
The author would like to thank Prof. Dr. Franz Pedit for his help, guidance and support in the research that produced this paper and Aaron Gerding for help in proof reading and editing.}

\bibliographystyle{unsrt}	
\bibliography{paper}

\begin{thebibliography}{10}

\bibitem{unf}
F.~Burstall and U.~Hertrich-Jeromin.
\newblock Harmonic maps in unfashionable geometries.
\newblock {\em Manuscripta Mathematica}, 108(2):171--189, 2002.

\bibitem{cec}
T.~Cecil.
\newblock {\em Lie Sphere Geometry}.
\newblock Springer, 2008.

\bibitem{bla}
W.~Blaschke.
\newblock {\em Vorlesungen {\"u}ber Differentialgeometrie {III}:
  Differentialgeometrie der Kreise und Kugeln}.
\newblock Springer, 1929.

\bibitem{sym}
A.~Sym.
\newblock Soliton surfaces and their applications (soliton geometry from
  spectral problems).
\newblock 239:154--231, 1985.

\bibitem{bur}
F.E. Burstall and U.~Hertrich-Jeromin.
\newblock The {R}ibaucour transformation in {L}ie sphere geometry.
\newblock {\em Differential Geometry and its Applications}, 24(5):503--520,
  2006.

\bibitem{bob1}
A.~I. Bobenko and Y.~B. Suris.
\newblock {\em Discrete Differential Geometry: Integrable Structure}, volume~98
  of {\em Graduate Studies in Mathematics}.
\newblock AMS, 2008.

\bibitem{dol6}
A.~Doliwa.
\newblock Discrete asymptotic nets and w-congruences in {P}l\"ucker line
  geometry.
\newblock {\em Journal of Geometry and Physics}, 39(1):9--29, 2001.

\bibitem{nov}
V.~E. Zakharov, N.~M. Ercolani, and F.~Calogero, editors.
\newblock {\em What is integrability?}
\newblock Springer Series in Nonlinear Dynamics. Springer, 1991.

\bibitem{dol2}
A.~Doliwa, P.M. Santini, and M.~Ma{\~n}as.
\newblock Transformations of quadrilateral lattices.
\newblock {\em J. Math. Phys.}, 41(944), 2000.

\bibitem{eis}
L.P. Eisenhart.
\newblock {\em Transformations of Surfaces}.
\newblock Oxford University Press, 1923.

\bibitem{dol3}
A.~Doliwa.
\newblock Quadratic reductions of quadrilateral lattices.
\newblock {\em Journal of Geometry and Physics}, 30(2):169--186, 1999.

\bibitem{yag}
I.M. Yaglom.
\newblock {\em Felix {K}lein and {S}ophus {L}ie: Evolution of the idea of
  symmetry in the nineteenth nentury}.
\newblock Birkhauser-Verlag, 1988.

\bibitem{dol1}
A.~Doliwa.
\newblock Generalized isothermic lattices.
\newblock {\em J. Phys. A: Math. Theor.}, 40(12539), 2007.

\bibitem{bob3}
A.~I. Bobenko, T.~Hoffmann, and B.A. Springborn.
\newblock Minimal surfaces from circle patterns: Geometry from combinatorics.
\newblock {\em Annals of Mathematics}, 164(1):231--264, 2003.

\bibitem{mob}
U.~Hertrich-Jeromin.
\newblock {\em Introduction to {M}{\"o}bius differential geometry}, volume 300
  of {\em London Mathematical Society Lecture Notes Series}.
\newblock 2003.

\bibitem{h}
F.~Burstall, D.~Ferus, K.~Leschke, F.~Pedit, and U.~Pinkall.
\newblock {\em Conformal Geometry of Surfaces in the 4-Sphere and Quaternions},
  volume 1772 of {\em Lectures Notes in Mathematics}.
\newblock Springer-Verlag, 2002.

\bibitem{ww}
R.S. Ward and R.~O.~Wells Jr.
\newblock {\em Twistor geometry and field theory}.
\newblock Cambridge Univ. Press, 1988.

\bibitem{bob2}
A.~I. Bobenko and Y.~B. Suris.
\newblock On organizing principles of discrete differential geometry: Geometry
  of spheres.
\newblock {\em Russian Math. Surveys}, 62:1--43, 2007.

\bibitem{udo2}
U.~Hertrich-Jeromin.
\newblock Transformations of discrete isothermic nets and discrete cmc-1
  surfaces in hyperbolic space.
\newblock {\em Manuscripta Mathematica}, 4:465--486, 2000.

\bibitem{sch}
F.~Burstall, F.~Pedit, and U.~Pinkall.
\newblock Schwarzian derivatives and flows of surfaces.
\newblock {\em Contemp. Math.}, 308:39--61, 2002.

\bibitem{pink}
U.~Hertrich-Jeromin, T.~Hoffmann, and U.~Pinkall.
\newblock A discrete version of the {D}arboux transform for isothermic
  surfaces.
\newblock {\em Inventiones Mathematicae}, 146(3):507--593, 1996.

\bibitem{bobpink}
A.~I. Bobenko and U.~Pinkall.
\newblock Discrete isothermic surfaces.
\newblock {\em Journal Fur Die Reine Und Angewandte Mathematik},
  47(5):187--208, 1996.

\bibitem{plu}
D.~Ferus, K.~Leschke, F.~Pedit, and U.~Pinkall.
\newblock Quaternionic holomorphic geometry: {P}l\"ucker formula, {D}irac
  eigenvalue estimates and energy estimates of harmonic 2-tori.
\newblock {\em Inventiones Mathematicae}, 146(3):507--593, 2001.

\bibitem{per}
U.~Hertrich-Jeromin, I.~McIntosh, P.~Norman, and F.~Pedit.
\newblock Periodic discrete conformal maps.
\newblock {\em Journal Fur Die Reine Und Angewandte Mathematik}, 534:129--154,
  2001.

\bibitem{dol7}
A.~Doliwa and P.M. Santini.
\newblock An elementary geometric characterization of the integrable motions of
  a curve.
\newblock {\em Physics Letters A}, 185(4):373 -- 384, 1994.

\bibitem{bohle}
C.~Bohle.
\newblock {\em M{\"o}bius invariant flows of tori in ${S}^4$}.
\newblock PhD thesis, T{\"U}-Berlin, 2003.

\bibitem{pp}
C.~Bohle, K.~Leschke, F.~Pedit, and U.~Pinkall.
\newblock Conformal maps from the 2-torus to the 4-sphere.
\newblock {\em in press}, 2011.

\bibitem{pen}
R.~Penrose.
\newblock {\em The Road to Reality}.
\newblock Johnathan Cape, 2003.

\bibitem{ped}
D.~Pedoe.
\newblock {\em Geometry, a comprehensive course}.
\newblock Dover publ., 1998.

\bibitem{bry}
R.~Bryant.
\newblock Surfaces in conformal geometry.
\newblock volume~47, pages 227--240, 1988.

\bibitem{cir}
H.~Schwerdtfeger.
\newblock {\em Geometry of Complex Numbers}.
\newblock Dover publ., 1980.

\bibitem{fer}
E.V. Ferapontov.
\newblock The analogue of {W}ilczynski's projective frame in {L}ie sphere
  geometry: {L}ie-applicable surfaces and commuting {S}chr{\"o}dinger operators
  with magnetic fields.
\newblock {\em International Journal of Mathematics}, 13(9):959--986, 2002.

\bibitem{dol4}
A.~Doliwa and P.M. Santini.
\newblock Multidimensional quadrilateral lattices are integrable.
\newblock {\em Physics Letters A}, 233(4-6):365 -- 372, 1997.

\bibitem{dol5}
A.~Cie\'{s}li\'{n}ski, A.~Doliwa, and P.M. Santini.
\newblock The integrable discrete analogues of orthogonal coordinate systems
  are multi-dimensional circular lattices.
\newblock {\em Physics Letters A}, 235(5):480--488, 1997.

\bibitem{bob4}
V.~E. Adler, A.~I. Bobenko, and Yuri~B. Suris.
\newblock Integrable discrete nets in grassmannians.
\newblock {\em Lett. Math. Phys.}, 89(2):131--139, 2009.

\bibitem{ks1}
B.~G. Konopelchenko and W.~K. Schief.
\newblock Menelaus' theorem, clifford configurations and inversive geometry of
  the schwarzian kp hierarchy.
\newblock {\em J. Phys. A: Math. Gen.}, 35(29), 2002.

\bibitem{ks2}
B.~G. Konopelchenko and W.~K. Schief.
\newblock Conformal geometry of the (discrete) schwarzian davey-stewartson ii
  hierarchy.
\newblock {\em Glasgow Mathematical Journal}, 47:121--131, 2005.

\bibitem{ks3}
B.~G. Konopelchenko and W.~K. Schief.
\newblock A novel generalization of clifford's classical point-circle
  configuration. geometric interpretation of the quaternionic discrete
  schwarzian kadomtsev-petviashvili equation.
\newblock {\em Proc. R. Soc. A}, 465(2104):1291--1308, 2009.

\end{thebibliography}


\end{document}